%% file: KROnset5.1.tex
\documentclass{eptcs}

\input{.//mydefs}

\input{.//journals}

\usepackage{amsmath, amssymb}
\usepackage{graphicx}
\usepackage{pstool}
\usepackage{mathtools}
\usepackage{psfrag}
\usepackage{hyperref}

\newcommand{\floor}[1]{\lfloor#1\rfloor}

\newcommand{\Nmin}{N_\mathrm{min}}

\newcommand{\ttherm}{T_\mathrm{therm}}
\newcommand{\tcorr}{\tau}
\newcommand{\cA}{\mathcal A}

\newcommand{\Integers}{{\hbox{\openface Z}}}

\begin{document}

\baselineskip = 20pt

\def\titlerunning{Asymptotic Regime for Finite Orders}
\def\authorrunning{J. Henson,  D. Rideout, R.D. Sorkin, S. Surya}

\title{Onset of the Asymptotic Regime \\ for Finite Orders}
\author{
Joe Henson$^\ast$, David Rideout$^\dagger$,
Rafael D. Sorkin$^\ddagger$
and Sumati Surya$^\circ$\\
$^\ast$ University of Bristol, U.K.  \\
$^\dagger$ University of California, San Diego, USA  \\
$^\ddagger$ Perimeter Institute, Waterloo, Canada \\
$^\circ$ Raman Research Institute, Bangalore, India
}
\date{\today}
\maketitle

\begin{abstract}
 We describe a Markov-Chain-Monte-Carlo algorithm which can be used to
 generate naturally labeled $n$-element posets at random with a
 probability distribution of one's choice.  Implementing this algorithm for the uniform
 distribution, we explore the approach to the asymptotic regime in which
 almost every poset takes on the three-layer structure described by Kleitman and
 Rothschild (KR).  By tracking the $n$-dependence of several
 order-invariants, among them the height of the poset, we observe an
 oscillatory behavior which is very unlike a monotonic approach to the KR
 regime.  Only around $n=40$ or so does this ``finite size dance''
 appear to give way to a gradual crossover to asymptopia which lasts
 until $n=85$, the largest $n$
 we have simulated.
\end{abstract}

\section{Introduction}

In the opposite regimes of small and large $n$,
much is known about the structure of the $n$-orders.
(By ``$n$-order'' we mean an $n$-element partial order or {\it\/poset},
or equivalently a finite $T_0$ topology on $n$ points.)
One the one hand, the unlabeled $n$-orders
(isomorphism equivalence classes of $n$-orders)
have been constructively enumerated by machine computation for $n\leq16$,
and their labeled counterparts
have been similarly enumerated to $n=18$ \cite{posets16}.
On the other hand, the limit of large $n$ has been treated by
Kleitman and Rothschild (KR) \cite{kr},
whose theorem tells us that as $n\rightarrow\infty$ the fraction of
posets belonging to a type which we will refer to
as {\it\/3-layered}
tends to unity.
This result yields rather complete information about
the structure of a typical $n$-order in the asymptotic regime,
including information on the sizes of the three layers,
as discussed below.
It implies in particular that
the total number of posets  of cardinality $n$
is to leading order
$2^{{n^2}/{4}}\,$,
independently of whether we consider the labeled or unlabeled case.

Intermediate values of $n$ are less well understood, however.  Not only
is it unknown how a typical poset is structured in this regime, but
it is not even known at which value of $n$ the asymptotic KR behaviour
sets in.
Of course no question like the latter can expect a precise answer.
Are we in asymptopia if 90\% of the $n$-orders are of KR type,
or must it be 99.9\%, or 99.99\%?
More importantly,
any attempt to define precisely the ``KR class'' of posets will
depend on which features one is interested in,
for example the height or the number of related pairs of elements.
Herein we consider a number of such indicators, and 
find a
relatively consistent picture emerging, of where and how the
transition to asymptopia takes place.

We address this question numerically, using the technique known as
Markov Chain Monte Carlo (MCMC).  Our MCMC algorithm is designed to
sample uniformly from the set $\Omega_n$ of naturally labeled
$n$-orders, where by a natural labeling we mean a labeling by natural
numbers which is compatible with the partial ordering itself.  (In
effect, the algorithm weights a poset by the number of its linear
extensions.  We have found such a weighting to be particularly easy to
implement numerically, and it arises naturally in the context of certain
growth dynamics for causal sets~\cite{dynamics,rideout}.)
While a random sampling is not as definitive as
an exhaustive enumeration,
it does provide
important
evidence on the structure of $\Omega_n$,
producing some surprising
(or at least, as far as we know, unanticipated)
results.

By design, our Markov dynamics satisfies ergodicity and detailed
balance, but that alone would only guarantee the desired sampling
probabilities in the unrealizable limit of an infinitely long run (and
provided we had a perfect random number generator).
To assess whether a uniform sampling has in fact been achieved, we have
run our Markov process with widely different initial posets, finding
excellent evidence that thermalization occurs on practical
time-scales
for $n\leq 85$.
Agreement with known exact results for $n=9$ provides a further reason
for confidence in our simulations.

A 3-layered poset has a height of three
(or less depending on the precise definition one adopts),
and height thus provides one possible indicator of when the asymptotic
regime has been reached.
Judged on this indicator,
our simulations exhibit a surprisingly colorful behavior
which lasts until $n\approx 45$,
and only then appears to switch over
to a gradual approach to the asymptotic limit.
To be more precise, the measured fraction of posets of height 3
actually \textit{decreases} for $7\leq n\leq30$,
falling below $3\%$ at $n=30$.
At this point, however, it
initiates a slow increase,
such that
the height 3 posets
make up over $90\%$ of the total for $n\geq80$.

Similar results are obtained for other indicators of KR-like behavior,
including the possession of a
layered structure,
the cardinalities of the levels
(to be defined below),
and
the total number of relations.
Taken as a whole, the
evidence from our simulations suggests a gradual crossover to
asymptopia
which begins around $n=43$, and
is still in progress at $n=85$,
the largest value we have simulated.
At what point one could
say with confidence that one had unequivocally entered the asymptotic regime is
difficult to assess from our data, but it looks to be at least $n=100$ or
greater.

The development of a reliable MCMC technique for posets is a uniquely
challenging problem that 
(in comparison with the Ising model for example) 
will probably require qualitatively new methods 
for its solution.
The need to maintain the very nonlocal constraint of transitivity  
presents novel difficulties because 
it introduces a highly ``inseparable'' mutual dependence 
among the relations defining the poset; 
and as far as we are aware
simulations of the uniform distribution on $\Omega_n$
have not been attempted until now
(although see \cite{mcmc2d} for a similar application of MCMC techniques).
On this first attempt, we have found our MCMC scheme to be
practical up to $n=85$, or for a perhaps fairer comparison with other
simulations, up to 3321 possible relations.  This offers hope that larger
values of $n$ could be accessed by further technical developments.

One of our main reasons for developing the techniques presented in this
paper comes from the poset's potential role in quantum gravity, and
specifically from its role in the theory of causal sets.
The point-events of a spacetime
which is free of causal pathologies\footnote{technically, a time-oriented Lorentzian geometry which
 contains no closed causal curve}
are partially ordered by
the relativistic light-cone structure.  This is a key motivation for
the causal set programme,
which
substitutes for
the spacetime continuum a locally finite
partially ordered set \cite{causets2}. 
Such a poset
thus represents a possible causal structure, and one can regard
a natural labeling thereof as a possible ``birth order'' of its
elements.
We expect the
techniques and
results presented herein
to be relevant to the study of
non-uniform measures on $\Omega_n$, such as arise \cite{BD_action} in
attempts to devise a quantum-gravity theory based upon the dynamics of
the discrete causal structures we have just described.

\section{Initial considerations}
\label{def.sec}
\subsection{Definitions}
\label{s:defs}

A partially ordered set or \emph{poset}
is determined by a pair $(S,\prec)$,
where $S$ (also called the \emph{ground set}) is a set
and $\prec$ is a binary relation on $S$ which is transitive
($x \prec y \prec z \implies x \prec z \;\forall x, y, z \in S$)
and (in our convention) irreflexive
($x \not\prec x \;\forall x \in S$).
When $x \prec y$,
we will sometimes express this by saying that
$x$ ``precedes'' $y$ or $y$ ``follows'' $x$.
We will consider
for
our simulations only
posets on the finite set $S = [n]$, where $[n] = \{0, 1, 2, \ldots n-1\}$,
and all of our posets 
will be
\emph{naturally labeled},
in the sense that
$x \prec y \implies x < y, \;\forall x, y \in [n] \;$.
We will refer this set of posets
as the \emph{sample space}
and denote it as $\Omega_n$.
An isomorphism equivalence class of posets will be called
an {\it\/unlabeled poset}.

By an \emph{order invariant} or \emph{observable} of a poset $S$ we will
refer to a property that is invariant under isomorphisms, or in other
words is independent of the labeling, natural or otherwise.  We will now
define a set of such invariants which will be employed in the
remainder of the paper.
An \emph{antichain} is a subset of $S$ among whose elements there are no relations $\prec\;$.
A \emph{chain} is a subset of $S$,
every pair of whose elements is
related by $\prec\;$,
and the {\it\/length} of a chain is its cardinality.
The \emph{height} of a poset is the length of the longest chain that it contains.
A \emph{minimal element} of $S$ is one with no element of $S$ preceding it,
while a \emph{maximal element} has no element following it.
A \emph{link} is a relation $x\prec y$
which is not implied by other relations via transitivity,
\emph{i.e.},
$x \prec y$ and $\nexists\; z$ such that $x \prec z \prec y$.
The total numbers of links $L$ and relations $R$ in a given poset
offer other useful invariants.
We define the \emph{linking fraction}
$l = 4L/n^2$ or $ 4L/(n^2-1)$ for $n$ even or odd, respectively,
and the \emph{ordering fraction} $r= R/\binom{n}{2}$,
whose denominator is the total number of possible relations on $n$ elements.
The \emph{past} of an
element $x$
consists of those elements which precede it in the order:
 $\past(x) = \{z : z \prec x\}$.
The \emph{inclusive past} of $x$ includes $x$ as well, $\ipast(x) = \past(x) \union \{x\}$.
Likewise we define the \emph{future} and \emph{inclusive future}
as $\fut(x) = \{z : x \prec z\}$
and $\ifut(x) = \fut(x) \union \{x\}$.
For $x\prec y$
the \emph{interval}
$I(x,y) = \fut(x) \cap \past(y)\,$.
Counting,
for each $k$,
the number of intervals of cardinality $k$
produces another set of invariants
containing a large amount of information about the structure of the
poset.

One can
recursively
partition a poset into \emph{levels} $L_1, L_2, \ldots, L_k$ as follows.
An element $x \in L_1$ if $\past(x)=\emptyset$, otherwise $x \in L_n$ if
the maximum value of $m$ for which there exists $y \prec x$
with $y \in L_m$ is $n-1\;$.
Equivalently, $n$ is the cardinality of the longest chain terminating at $x$.

\subsection{The Kleitman-Rothschild theorem}
\label{s:krtheorem}

The Kleitman-Rothschild theorem offers a precise answer to
the question of
what a typical poset looks like
in the limit of large $n$.
We state it in essentially the same form as Brightwell \cite{brightwell},
making use of the following
families  of ``layered'' posets.
Consider
disjoint sets $X_1, X_2, \ldots, X_k$ of elements, and let
$\mathcal{A}(X_1, X_2, \ldots, X_k)$ denote the set of partial orders
definable on $\union_{i=1}^k X_i$ such that
\begin{enumerate}
  \item[(1)] if $x\in X_i$, $y\in X_j$, and $x \prec y$, then $i < j$, and
  \item[(2)] if $x\in X_i$, $y\in X_j$, and $i < j-1$, then $x \prec y \;$.
\end{enumerate}
Note that each $X_i$ forms an antichain, which we call \emph{layer} $i$.
For any poset whatsoever,
its partition into {\it\/levels} clearly fulfills condition (1),
but not necessarily (2).
Finite posets in the class $\mathcal{A}(X_1, X_2, \ldots, X_k)$ have a number of
readily deduced properties,
the most obvious of which is that they have height $ \leq k$.

We will say that \textit{in the limit of large $n$, almost every poset}
has a certain property if, as $n \rightarrow \infty$, the fraction of
posets lacking that property goes to $0$.
We will also say that the posets with such a property \textit{dominate at large $n$}.

\paragraph{KR theorem:} Let $\omega(n)$ be any function tending to infinity.
 In the limit of large $n$, almost every partial order with ground-set $[n]$ lies
 in the class $\cA(X_1,X_2,X_3)$ for some partition $(X_1,X_2,X_3)$ of $[n]$
 such that
 $||X_2|-n/2| < \omega(n)$,
 $||X_1|-n/4| < \omega(n)\sqrt{n}$, and
 $||X_3|-n/4| < \omega(n)\sqrt{n}$.

\paragraph{}
The above theorem informs us that asymptotically almost every $n$-order
is a \emph{three-layer poset} in the sense that it belongs to the class
$\cA(X_1,X_2,X_3)$, and from this follow other asymptotic properties not
stated explicitly in the theorem.
Such properties are often easy to derive, if one notices that in the
uniform random distribution of posets in $\cA(X_1,X_2,X_3)$, each random
variable associated to the existence of a relation between elements in
$X_1$ and $X_2$, or between $X_2$ and $X_3$, is independent.
For example, one can prove in this way that
the posets of height 3 dominate at large $n$.
This is not true tautologically,
because our definition of a layered order allows any of the three layers to be empty,
and because condition (2) above does not force
an
element in layer 2 to be related
to any other element in layer 1 or 3.
Nonetheless, in view of the previous observation, it is not difficult to
see that the fraction of 3-layer posets of height $< 3$
tends to 0 as $n\rightarrow\infty$.
By similar reasoning
it is not hard to demonstrate that,
for the uniform distribution over 3-layer
orders, the ordering fraction
is concentrated near to $r=3/8$.
From the theorem as stated, one can also derive,
by counting the possible 3-layer orders,
that the total number of $n$-orders grows as $2^{n^2/4}$,
up to a subdominant supplementary factor.
All of these properties
can be regarded as characteristic of the asymptotic regime.

The $n$-orders referred to in the KR theorem as stated above
are ``arbitrarily labeled posets'',
\textit{i.e.}~posets on $[n]$ with no further restrictions.
Although the
set of all such posets
will obviously differ from the set of
unlabeled $n$-orders,
Kleitman and Rothschild showed that
their result holds for unlabeled orders as well.
Intuitively,
a labeling can introduce a factor of at most $n!$,
and $n!$,
imposing as it might seem,
is subdominant compared to $2^{n^2/4}$.
(In this connection, it's worth emphasizing that for large $n$, almost every
3-layer $n$-order is asymmetric: it admits no nontrivial automorphism.)

As mentioned, however, our simulations produce neither unlabeled nor
arbitrarily labeled posets, but naturally labeled ones, a choice which
allows for efficient machine representation of a poset as an
upper-triangular matrix.\footnote{An upper-triangular bit-matrix $M$ defines a binary relation on $[n]$
  via the rule $j\prec k\iff M_{jk}=1$.  If one adds the condition of
  transitivity, then $M$ yields a general element of $\Omega_n$.}
In comparison with a uniform weighting over unlabeled $n$-orders, our
sample-space thus weights each unlabeled poset by the number of its
linear extensions (or more precisely by the number of inequivalent
linear extensions, where two linear extensions are equivalent iff they
are related by an automorphism of the poset).

Insofar as one cares only about
enumerating the $n$-orders asymptotically,
the distinctions among unlabeled, naturally labeled, and
arbitrarily labeled posets can (as just explained) be ignored, unless
one is interested in the subleading terms:
the resulting under- or
over-counting could ``at most'' 
introduce
a factor of $n!$.  
However, 
the KR theorem,
in the form given above, 
also says something about the typical sizes of the
three layers, and here the choice of natural labeling will make a
difference, because a reallocation of poset elements between the top
and bottom layers has no effect on the leading order counting.
A closer analysis,
based on an estimate we reproduce in the Appendix,
leads us to expect that
the size of the bottom (or top) layer is more likely to vary uniformly
between 0 and $n/2$ than to be concentrated near $n/4$ as the KR theorem
would have it.  And this in turn would imply a less concentrated
asymptotic ordering-fraction averaging to $1/3$ instead of $3/8$.
(Asymptotically $r$ would be distributed between 1/4 and 3/8 with a
square-root-divergent peak at 3/8, which accordingly would still be its
``most likely value''.)

Finally, let us comment that the choice of natural labeling can also be
motivated by its application to the physics of spacetime, and
specifically by the example of ``sequential growth dynamics'' in causal
set theory.  In that context a measure is defined on the space of
countable posets corresponding to a certain Markov process, wherein the
labels acquire a temporal meaning \cite{dynamics, rideout, prob_of_time, observables}.

\subsection{In search of the asymptotic regime}
\label{s:onset}

As
stipulated
earlier, the question of where the asymptotic regime
begins is inherently ambiguous, but if we don't demand too much
quantitative precision, we can come up with some reasonably natural
criteria.
The question amounts to deciding what should count as a ``KR order'',
and then asking at what value of $n$ the majority of $n$-orders are in
fact of this form.

The theorem quoted above suggests first of all that one should define a
KR order to be 3-layered (albeit even this is subject to some doubt, as
the definition given above is only one of several possible variations on
the same theme).
The theorem also limits the sizes of the three layers, but unfortunately
it does so in a manner that tends to lose its meaning when applied to a
definite, finite $n$.  (Whether or not any given poset satisfies
$||X_2|-n/2| < \omega(n)$ depends entirely on what choice we make of the
arbitrary function $\omega(n)$.)
What remains is the semi-quantitative criterion that the middle layer
should hold about half the elements, with the other two layers each
holding around a fourth (albeit with larger fluctuations than for the
middle layer).  We have seen, however, that the condition on the top and
bottom layers is inappropriate when the sample space consists of
naturally labeled orders, as it does for us.  Instead we expect
something more like a uniform distribution of sizes for those layers.
Further natural criteria we have encountered are that the height should be
exactly 3 and that the ordering fraction $r$ should take on
a characteristic value, plausibly varying between 1/4 and 3/8 with a
mean of 1/3.
In the following, we explore all these criteria.

\section{The Markov chain}

The canonical application of the MCMC technique concerns a collection of
$\Integers_2$-valued ``spin'' variables on a regular lattice,
where the probability distribution one is attempting to simulate is
given by an easily computed ``Boltzmann weight''.
Two features that make such models easy to deal with are
that
(\textit{a}) the configurations are expressible in terms of a fixed
set of independent variables simply related to the Boltzmann weight,
and
(\textit{b}) there are available obvious ``local Monte Carlo moves'' that need only
refer to an easily defined, small region of the overall configuration,
\textit{e.g.} ``flipping'' one of the spins.
Not all successful MCMC simulations have relied on property (\textit{a}),
for instance one has simulated ensembles
of simplicial manifolds (see \textit{e.g.} \cite{cdts})
but in the case of posets even property (\textit{b}) is lacking.
As a result, efficient moves seem to become increasingly hard to
devise as $n$ increases.

In \cite{mcmc2d} MCMC methods were used to study the space
$\Omega_n^{2d}$ of $n$-element two dimensional orders, which is a proper
subset of $\Omega_n$.  The asymptotic limit for the uniform distribution
on 2d orders was found by Winkler and El Sauer
\cite{zaharsauer,winkler}, and an analogous question about onset of the
asymptotic regime could be asked in that case.
The results of \cite{mcmc2d}
indicate that the asymptotic regime is realized for $n$ as small as
$30$.  However, 2d orders are structurally very different from generic
posets, and it would be no surprise
if
the onset of asymptotic behavior
took on a different character in the two cases.
Notice also that features (\textit{a}) and (\textit{b}) are both present
for 2d orders, so that one learns
from them
little
about the process of equilibration
which could
be applied to the full $\Omega_n$.

The Markov chain on $\Omega_n$ which we used in the present work employs
a mixture of {\it\/relation moves} and {\it\/link moves}, as described in detail
below.
Roughly, one
removes or adds a relation or a link (as defined above) between a pair of
randomly chosen elements of the poset.
The simulations reported herein used a uniform mixture
of the two types of moves.
Each type was tried in independent trials
and it was found that thermalization was substantially
hastened
when
both types were  employed together.\footnote
{We  also experimented with a
  third move inspired by the dynamics of ``transitive percolation'' \cite{dynamics}.
  However it did not seem to speed up the thermalization of the Markov chain.}

In order that a Markov chain on some sample-space of ``states'' converge to the desired equilibrium distribution,
it suffices that
the moves it employs be ergodic and that their conditional probabilities
be chosen to satisfy detailed
balance~\cite{metropolismc}.
Ergodicity merely requires that it be possible to pass (with non-zero
probability) from any state to any other state.
Detailed balance is the condition that
\bne
      \Pr(A) \Pr(A \to B) = \Pr(B) \Pr(B \to A) \;, \label{detbal.eqn}
\ene
where $A$ and $B$ are any two states,
$\Pr(X)$ is the desired probability of state $X$,
and $\Pr(X \to Y)$ is the transition-probability from $X$ to $Y$
that one specifies in setting up the Markov chain.
Since in our case we seek the uniform distribution over states
(i.e. posets in $\Omega_n$), we need only ensure that each move
has the same transition probability as
its inverse.
(We remark here that while ergodicity is trivially necessary for the
desired equilibrium, detailed balance is only a matter of choice.  Other
conditions would serve as well, and in certain situations they might be
more efficient.)

By a uniform mixture of relation- and link-moves
we simply mean that at each step we flip a fair coin
and propose one or the other type of move with equal probability.
If the transition
probabilities for the two move-types are given by $\Pr_r\;$ and $\Pr_l\;$,
then detailed balance for an equal mixture of the two types plainly
requires that
\bne
  \mathrm{Pr}_r(A \to B) + \mathrm{Pr}_l(A \to B) = \mathrm{Pr}_r(B \to A) +
  \mathrm{Pr}_l(B \to A) \;,
  \label{mixture.eqn}
\ene
a condition which is clearly satisfied if each move individually
satisfies its own detailed
balance condition (\ref{detbal.eqn}).
Likewise, a mixture of two ergodic move types is a fortiori ergodic.

In order to describe the relation and link moves fully we will, for the
remainder of this section, think of the elements of our $n$-order as
being the natural numbers from $0$ to $n-1$.  (Thus each element serves
as its own label.)
On the computer
we represent a labeled partial order by
its so-called adjacency matrix $A$, defined by:
$A_{ij} = 1$ iff $i \prec j$, otherwise 0.
Irreflexivity then sets the diagonal of $A$ to zero,
while the labeling is natural iff $A$ is upper triangular.

We will need two more definitions in order to describe our moves.
A \emph{critical pair} is a pair of elements $x<y, x \nprec y$
for which
$\past(x) \subseteq \past(y)$ and $\fut(y) \subseteq \fut(x)$.
A \emph{suitable} pair is a pair of elements $x<y, x \nprec y$,
such that there exists no $z\in\ipast(x)$
which is linked to an element $w\in\ifut(y)$.\\

The {\bf \emph{relation move}}  is:
\begin{enumerate}
  \item Select a pair of elements $x<y$ uniformly at random.
  \item If $x$ and $y$ are linked, remove the single relation $x \prec y$.
  \item Else if $x$ and $y$ form a critical pair, adjoin the single relation $x \prec y$.
  \item Otherwise leave the poset unchanged.
\end{enumerate}
Since a succession of relation-moves can clearly transform any poset to (say) a chain, they can
transform any poset to any other, hence this move is ergodic.
Moreover, it satisfies
detailed balance:
Consider first a linked pair $(x, y)$ whose relation $x \prec y$ gets
removed.
The probability to select the pair is $\frac{2}{n(n-1)}$,
while the probability to accept the move is 1.
After the move is completed $(x,y)$ will be a critical pair.
The inverse of this move begins with a critical pair $(x, y)$ and adjoins
the relation $x \prec y$.
The pair  $(x,y)$
is selected with the same probability $\frac{2}{n(n-1)}$,
and the move is accepted with probability 1.
After the move is completed $(x,y)$ will be a linked pair.
In the event that the selected pair forms neither
a link nor a critical pair, the poset does not change.
In this case the move is its own inverse, and ``both'' trivially share
the same probability.\\

The {\bf \emph{link move}} is:
\begin{enumerate}
 \item Select a pair of elements $x<y$ uniformly at random.
 \item If $x$ and $y$ are linked
 \begin{enumerate}
   \item Remove all relations between elements of $\ipast(x)$ and elements of $\ifut(y)$.
   \item Restore every relation which is implied by transitivity (via an element
         which is unrelated to at least one of $x$ or $y$.)
  \end{enumerate}
  \item Else if $x$ and $y$ form a suitable pair,
        adjoin a relation $x'\prec y'$ for every pair of elements $x'\in\ipast(x)$ and $y'\in\ifut(y)$.
  \item Otherwise leave the poset unchanged.
\end{enumerate}
If this description seems a bit complicated, it is because a link move acts directly,
not on the adjacency matrix but on the ``link matrix''.  That is it
effectively acts on the poset construed as a set of
{\it\/links} rather than as a set of {\it\/relations}.
Pictorially, it inserts or deletes an edge in the Hasse diagram of the poset,
which is a directed graph containing
an edge for every linked pair of poset-elements.
Since such diagrams are in
one to one correspondence with partial orders, it is clear that the link
move is ergodic.
(One can convert any poset to an antichain by removing all of its links.)
Detailed
balance is satisfied much in the same way as for the relation move.
The selection probability is $\frac{2}{n(n-1)}$ for each potential move
from a given poset.
Each potential move either leaves
the poset unchanged, or converts it to a different poset with
probability 1.
In the latter case,
the inverse move, which is the only way to return to the original poset,
has the same selection
and acceptance probability as the original move.
Thus detailed balance is satisfied.

Note that a relation added between a critical pair creates no further new
relations, and likewise a link inserted between a suitable pair creates no
further links.  (In fact, we could have taken this as the definitions of
such pairs.)
For this reason, the matching inverse moves are easy to identify as
simply deleting the relation or link that was added.

\section{Code details}

We have implemented the above Markov chain Monte Carlo algorithm as a
module within the \texttt{Cactus} high performance computing framework
\cite{cactus}, making use of the \texttt{CausalSets} toolkit (some details
are available in \cite{ccl}).  The \texttt{CausalSets} toolkit provides basic
infrastructure for working with partial orders, including data structures and
numerous algorithms needed to implement the above moves and compute order
invariant `observables'.

We compute pseudo-random numbers using the maximally equidistributed combined
Tausworthe generator of L'Ecuyer~\cite{gsl_rng_taus2}, as implemented in the
GNU Scientific Library~\cite{gsl}. 

For a given value of $n$, we generally start the Markov chain with four different starting posets,
to aid in determining when it thermalizes.  The starting posets are a chain, antichain, a random
Kleitman-Rothschild order (with $\floor{n/2}$ elements in the middle layer, the cardinality of the
bottom layer selected from a Poisson distribution with mean $\floor{n/4}$, and the remaining
elements going into the top layer), and
a bipartite order from $\mathcal{A}(X_1, X_2)$ with $|X_1| = \floor{n/2}$ and every
element of $X_1$ related to every element of $X_2$.
We then execute a large number of sweeps
(the exact number depending upon $n$),
where each sweep consists of $2n^3$ attempted moves.
At the
end of each sweep we record a large number of observables including height, ordering fraction, and numbers
of minimal and maximal elements.
We regard attempted moves which change the poset as `accepted', and others as `rejected'.  Each move
has an acceptance rate close to 1/2 for all the simulations that we have performed.
Note that 
 our choice of
 $2n^3$ attempted moves per sweep is longer than the more
 natural choice of ${n \choose 2}$ attempted moves, the number needed so
 that on average each relation is visited once per sweep.  
 We found,
 however,
 that sampling 
 (and then storing)
 the observables  
 as often as this 
 consumed a great deal of disk space, making sweeps of length $2n^3$ more
 practical. 
 (In fact, for $n\geq 69$, the autocorrelation times are long enough that 
 sweeps of length $2n^5$ 
 are justified.)

The results come from many core-years of compute time, both on the
Lonestar cluster at the Texas Advanced Computing Center (TACC), and a 12
core Intel Xeon X5690 workstation. In addition, extensive preliminary
investigations with the link move were done on the HPC cluster at the
Raman Research Institute.

Each realization of the Markov chain begins with a transient portion
which must be discarded because the probabilities have not yet reached
their equilibrium values.  They have not yet ``thermalized'', as one
says.  In deciding how long to wait, we employ as indicators the number
of minimal elements $\Nmin$ and the ordering fraction $r$, as these two
observables seem to be particularly sensitive to the approach to
equilibrium.  We regard the Markov chain as having thermalized when the
traces of these observables versus sweep behave similarly from all four
starting posets.  Examples of traces of $\Nmin$ and $r$ vs.\ sweep are
depicted in Figures \ref{obs_trace47.fig} and \ref{obs_trace67.fig}.
\begin{figure}[htb!]
  \begin{minipage}[t]{.5\textwidth}
    \includegraphics[width=\textwidth]{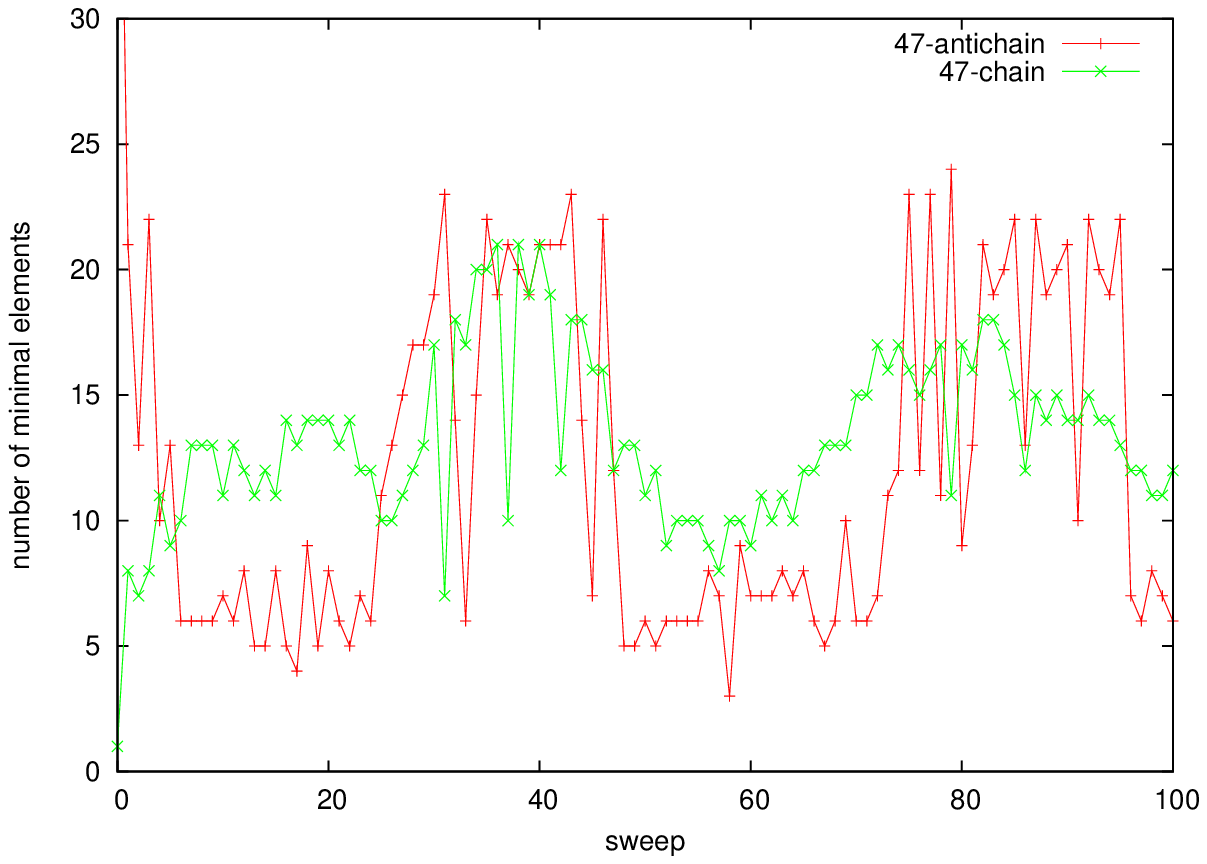}
  \end{minipage}
  \begin{minipage}[t]{.5\textwidth}
    \includegraphics[width=\textwidth]{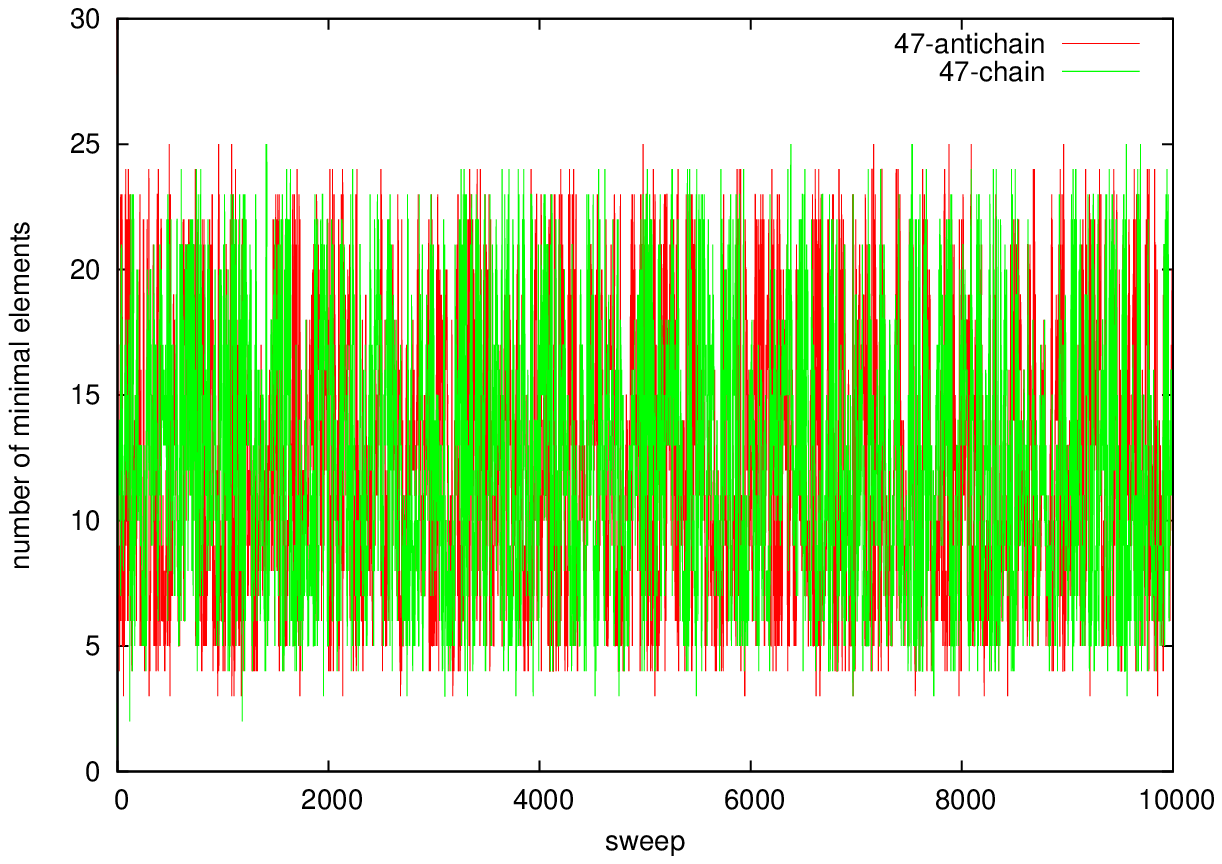}
  \end{minipage}
  \caption{\label{obs_trace47.fig} Trace of number of minimal elements for
     $n=47$, starting from an antichain (red) and a chain (green).  The Markov
     chain appears to be approaching equilibrium after 30 sweeps (if not sooner).}
\end{figure}
\begin{figure}[htb!]
  \begin{minipage}[t]{.5\textwidth}
  \includegraphics[width=\textwidth]{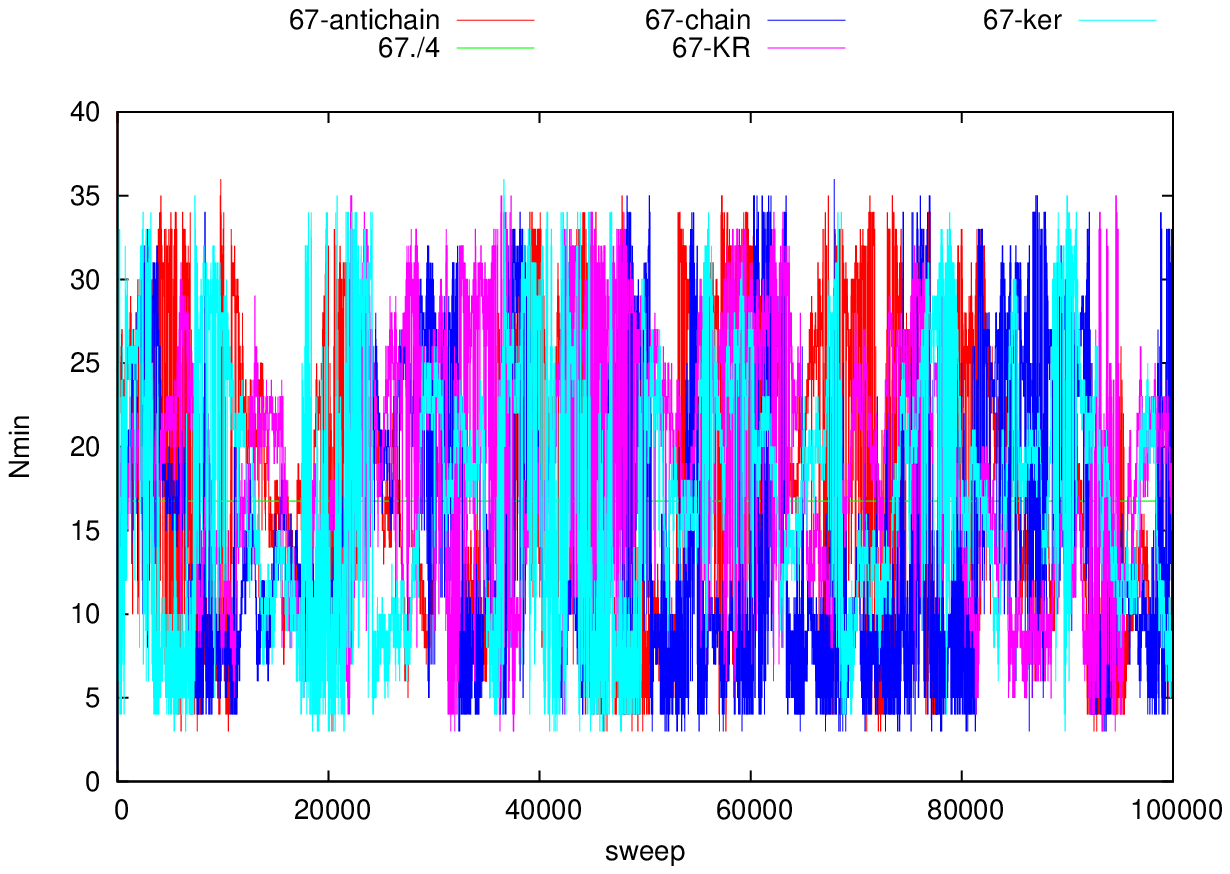}
  \end{minipage}
  \begin{minipage}[t]{.5\textwidth}
  \includegraphics[width=\textwidth]{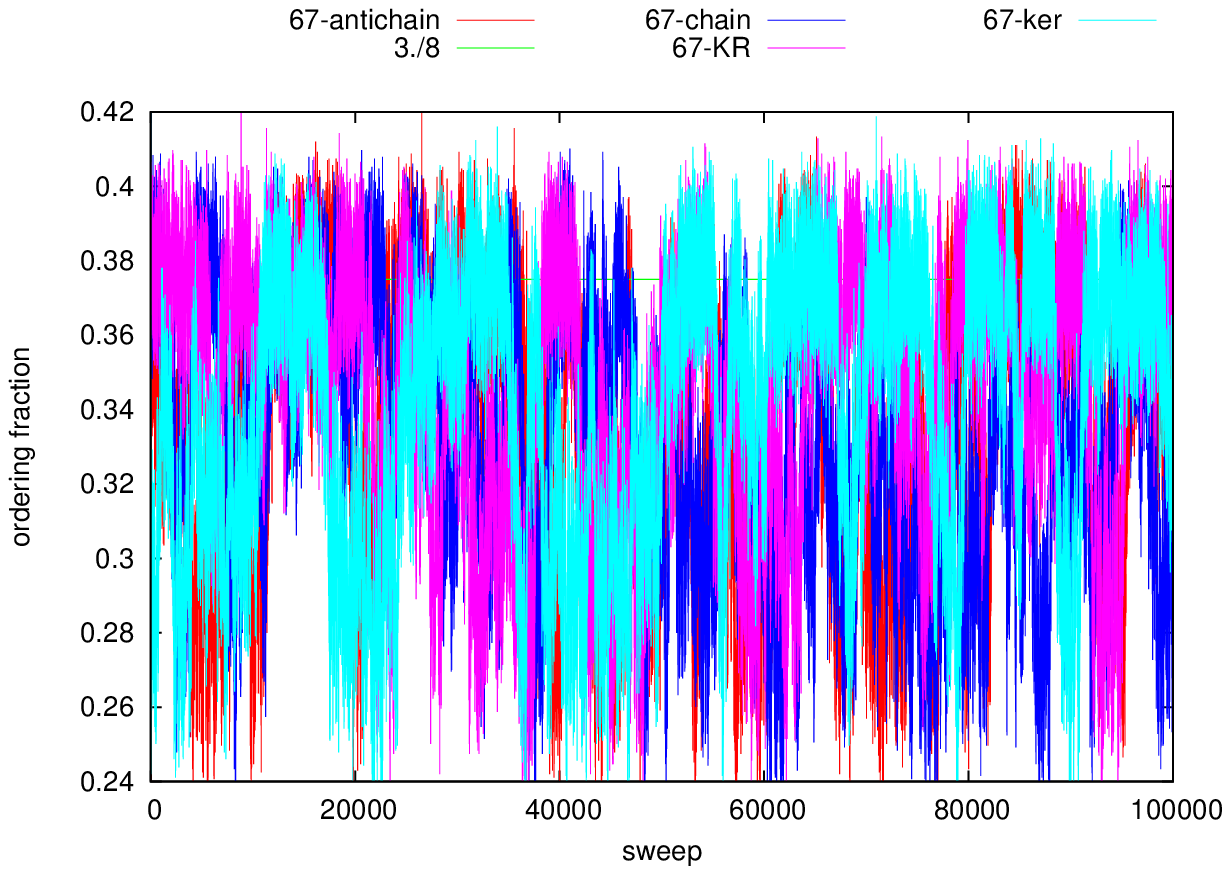}
  \end{minipage}
  \caption{\label{obs_trace67.fig} Trace of number of minimal elements (left)
    and the ordering fraction (right) for
    $n=67$, starting from an antichain (red), a chain (blue), a
    Kleitman-Rothschild order (magenta), and a bipartite order with linking
    fraction 1 (cyan).  The expected number of minimal elements for a Kleitman-Rothschild
    poset is $n/4$.
    The Markov chain appears to have thermalized after 41,000 sweeps.}
\end{figure}
These two observables reach equilibrium more slowly than 
others, 
such as 
the linking fraction, 
which appears thermalized 
even when the Markov chain is far from equilibrium.
This can be observed by comparing with, for example,
the number of minimal elements, as demonstrated in Figure \ref{linkfrac.fig}.
\begin{figure}[htb!]
\begin{minipage}[t]{.5\textwidth}
\includegraphics[width=\textwidth]{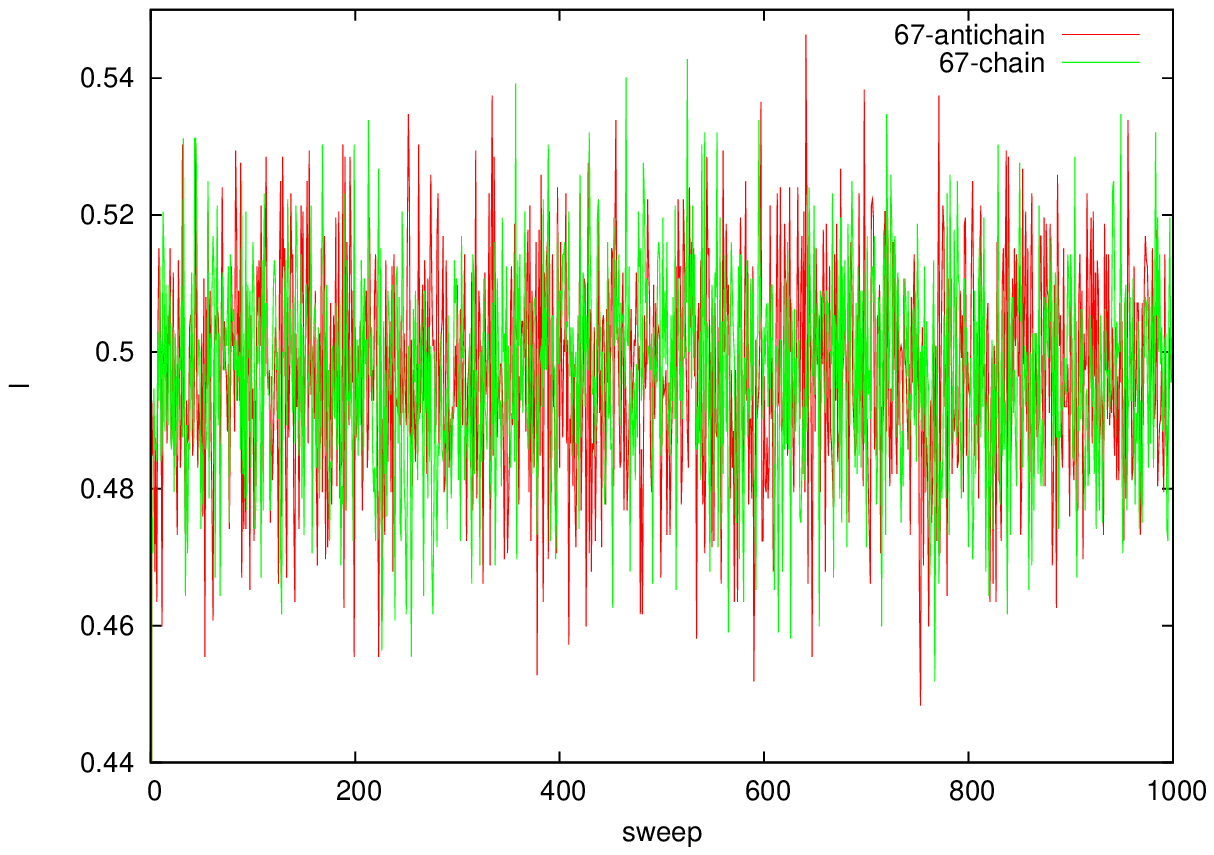}
\end{minipage}
\begin{minipage}[t]{.5\textwidth}
\includegraphics[width=\textwidth]{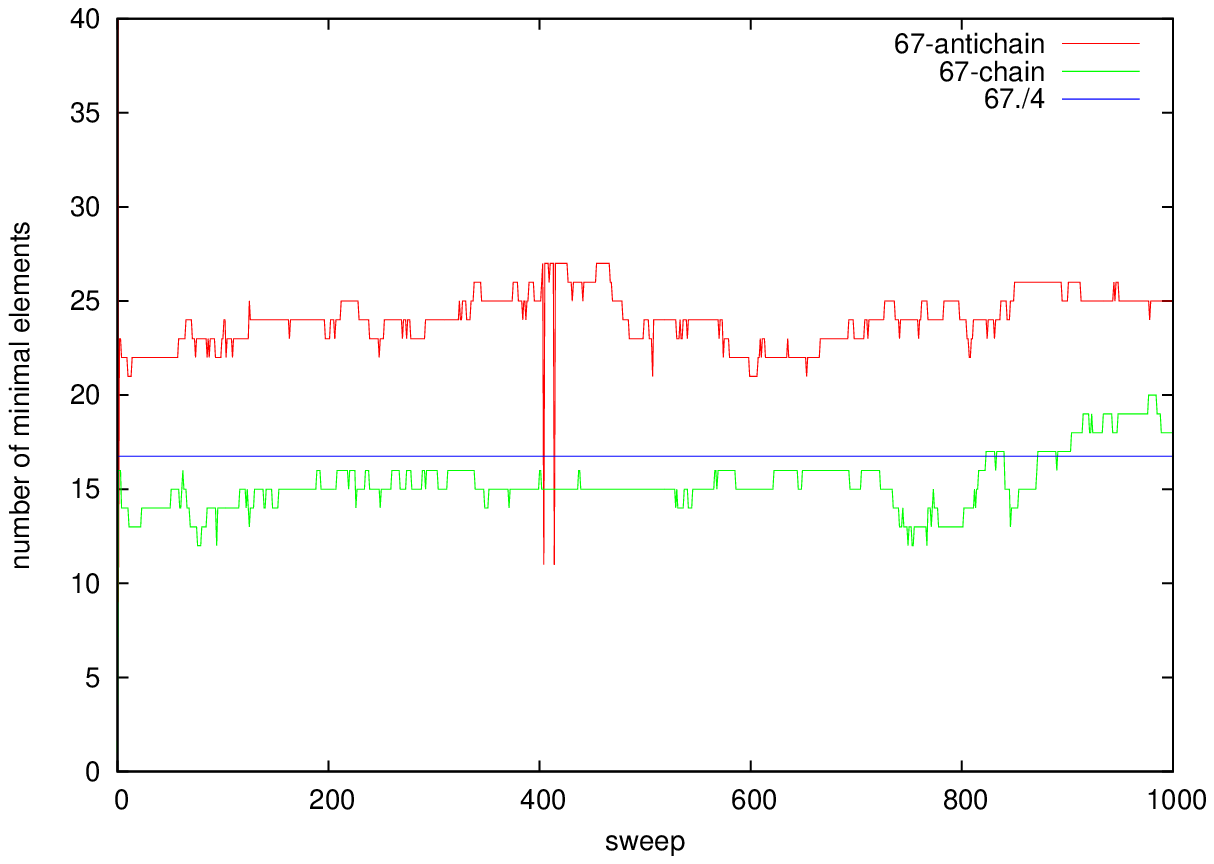}
\end{minipage}
\caption{\label{linkfrac.fig} Trace of linking fraction (left) and number of
  minimal elements (right) for the first 1000 sweeps at $n=67$,
  starting from an antichain (red) and a chain (green).
  By way of comparison, the thermalization time deduced from our other
  indicators is actually $\ttherm = 47,000$ for $n=67$.
  The expected linking fraction for a Kleitman-Rothschild poset is 1/2.}
\end{figure}

After discarding the initial, unthermalized region, 
we estimate the autocorrelation time $\tcorr$ in equilibrium,
by fitting,
for our various observables $O$,
the correlator,
$\langle (O(s) - \langle O \rangle) (O(s+t) - \langle O \rangle) \rangle$ 
(where $\langle O \rangle$ is the value of $O$ averaged over all the measurements),
as a function of lag-time $t$
against an exponential curve, $\;a\,\exp(-t/\tcorr)$.
Sample fits for $n=67$ are shown in Figure \ref{tau67.fig}.
\begin{figure}[htb!]
\psfrag{Nmin}{$\Nmin$}
\psfrag{x}{$t$}
\centerline{\includegraphics{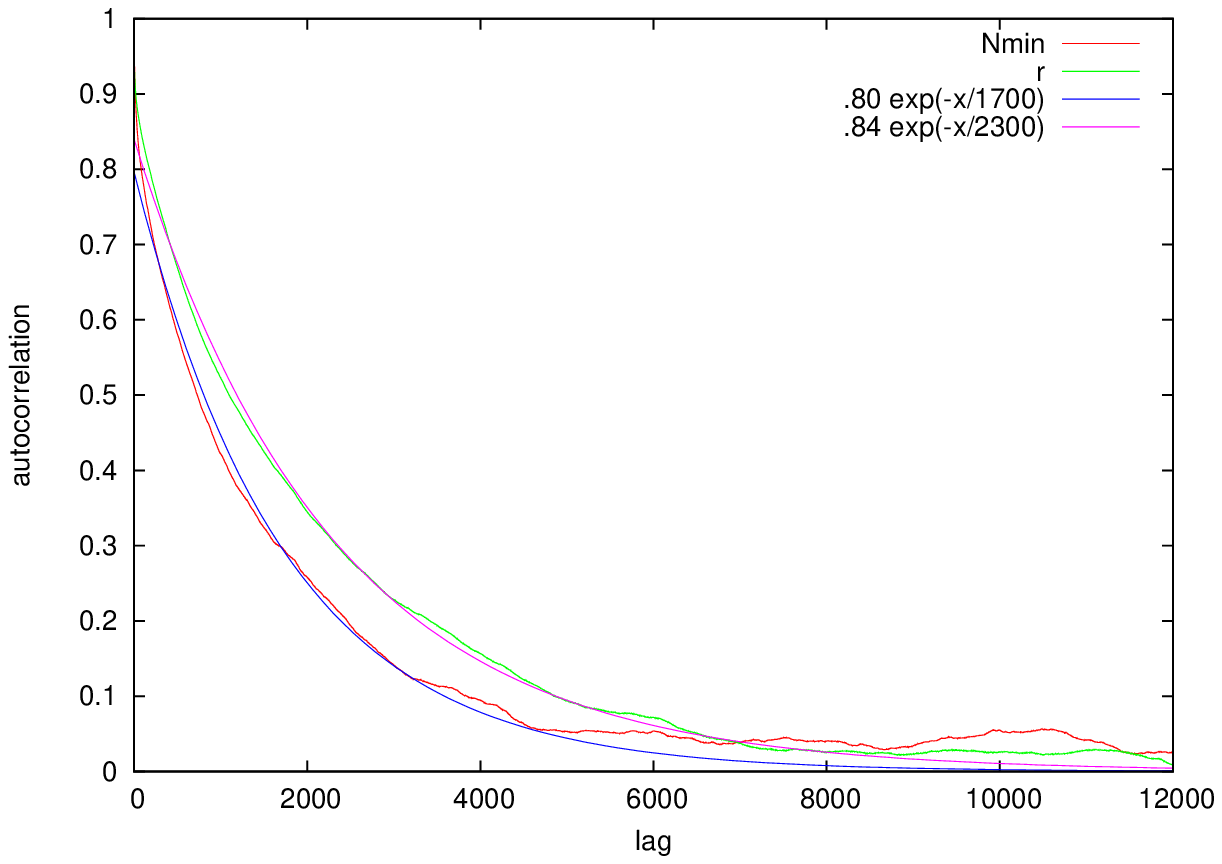}}
\caption{\label{tau67.fig} Autocorrelation functions of
  number of minimal elements (red) and ordering fraction (green) for $n=67$,
  fitted with exponential curves in order to measure the autocorrelation
  time in equilibrium.
  Here we used $\ttherm = 47,000$.}
\end{figure}

Figure \ref{acorr.fig} plots our estimates for thermalization time $\ttherm$
and autocorrelation time $\tcorr$, measured in sweeps.
Each appears to grow exponentially with $n$.
We performed this analysis for a limited set $\Phi$ of poset sizes $n$,
as illustrated in the figure.  In our simulations, however, we needed to
know the thermalization time $\ttherm$ for every $n$ that we simulated,
not only for those $n$ in $\Phi$.  For simplicity in such cases, we
estimated $\ttherm(n)$ conservatively from $\ttherm(n_0)$ for the
smallest $n_0\in\Phi$ such that $n_0 \geq n$,  
given the natural assumption that  thermalization times increase with $n$. 
\begin{figure}
\psfrag{Ttherm}{\small $T_\mathrm{therm}$}
\psfrag{Tcorr}{\small $T_\mathrm{corr}$}
\centerline{\includegraphics{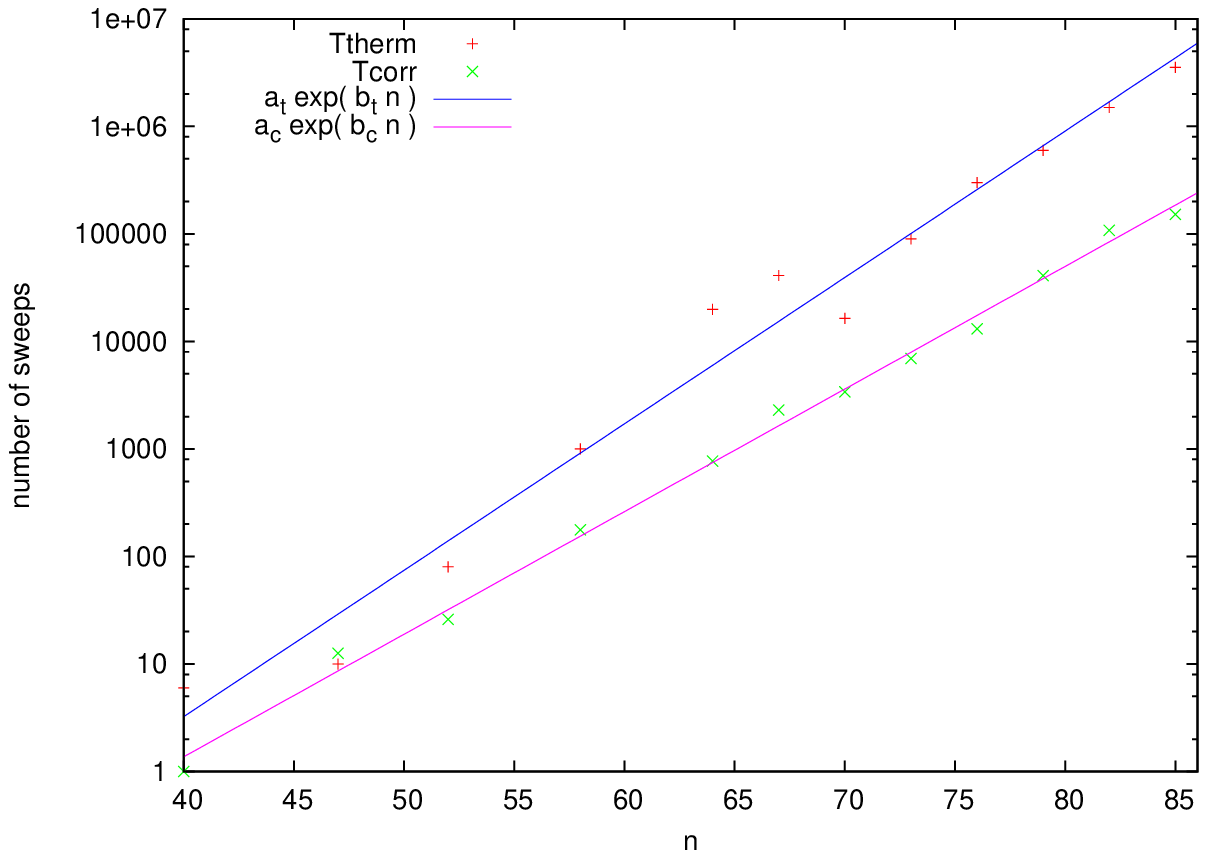}}
\caption{\label{acorr.fig} Measured values of thermalization time $T_\mathrm{therm}$ and
  autocorrelation time $T_\mathrm{corr}$ for various poset sizes $n$.  For $n<40$ the
  autocorrelation time appears to be less than a single sweep.
  The best fitting exponential functions,
  $\ttherm = a_t \exp(b_t n)$ and $\tcorr = a_c \exp(b_c n)$,
  have $\ln a_t = -11.4 \pm 1.0$, $b_t = 0.314 \pm 0.015$,
  $\ln a_c = -10.2 \pm .4$, and $b_c = 0.263 \pm .005$.}
\end{figure}

Equipped with a thermalization time $\ttherm$ and autocorrelation time $\tcorr$,
we compute histograms of a number of
order-invariants
such as height and number of relations, as follows.
The simulation outputs a sequence of observable
values after each sweep.
We discard $\ttherm$ sweeps, and use
the remaining samples
(taken at every $5 \tcorr / 2$
sweeps for $n>40$, and at every sweep otherwise)
to calculate the mean and its error.
Much of our subsequent analysis is based on
histograms built from these samples.
We estimate each bin frequency $f$ by the number
of events which fall into that bin divided by the total number of samples $T$,
and we estimate the statistical error in this frequency as $\sqrt{f(1-f)/(T-1)}$.
Experimentation indicates that the resulting error estimates adequately
reflect the uncertainty in the measured histogram bin frequencies.
Alternatively, one could sample the time series more frequently, and use
either the bootstrap or jackknife method to estimate the bin frequencies
and their errors

\section{Correctness}

One of us (RDS)
has generated a library of the unlabeled posets of nine or fewer
elements, and has counted the natural labellings for each of them
(details will appear separately \cite{smalln}).
We use this collection to check that our code is giving correct results.

\begin{figure}[htb!]
  \begin{minipage}[t]{.5\textwidth}
  \includegraphics[width=\textwidth]{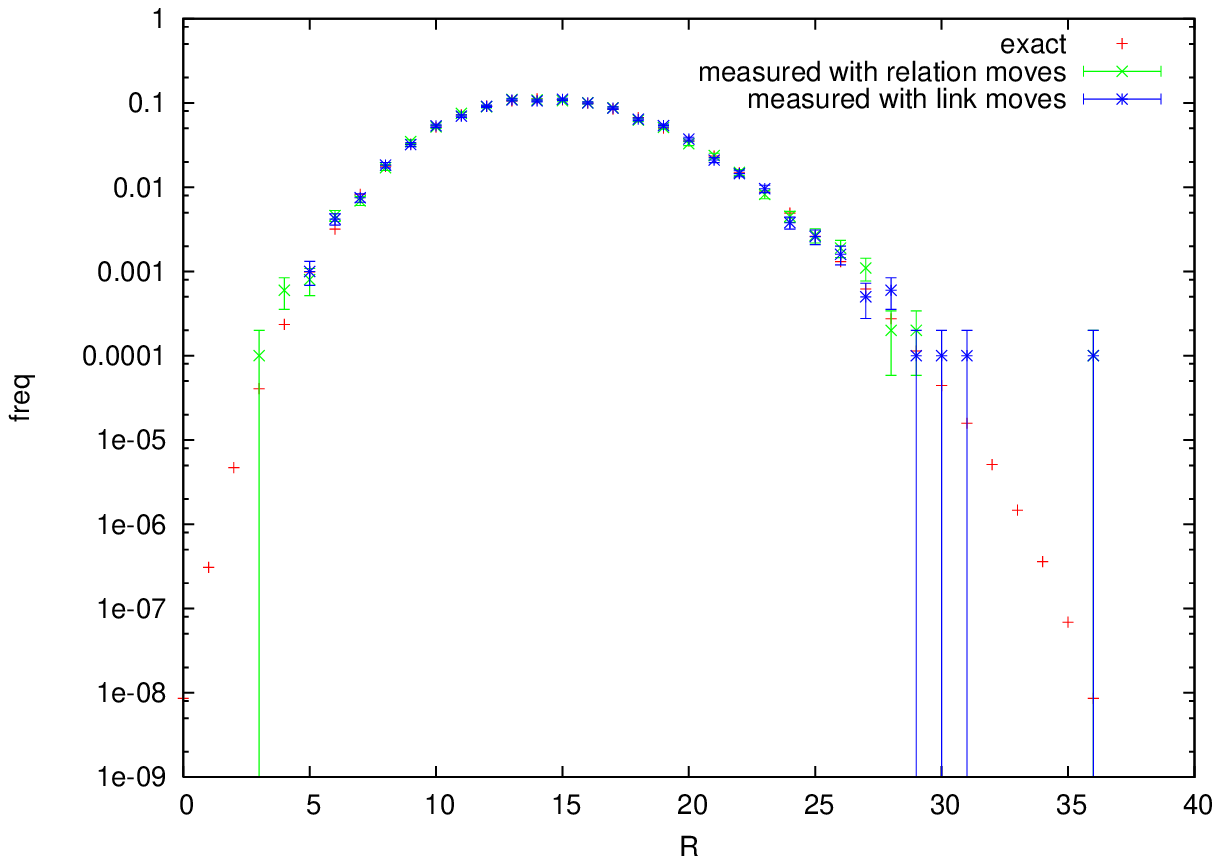}
  \end{minipage}
  \begin{minipage}[t]{.5\textwidth}
  \includegraphics[width=\textwidth]{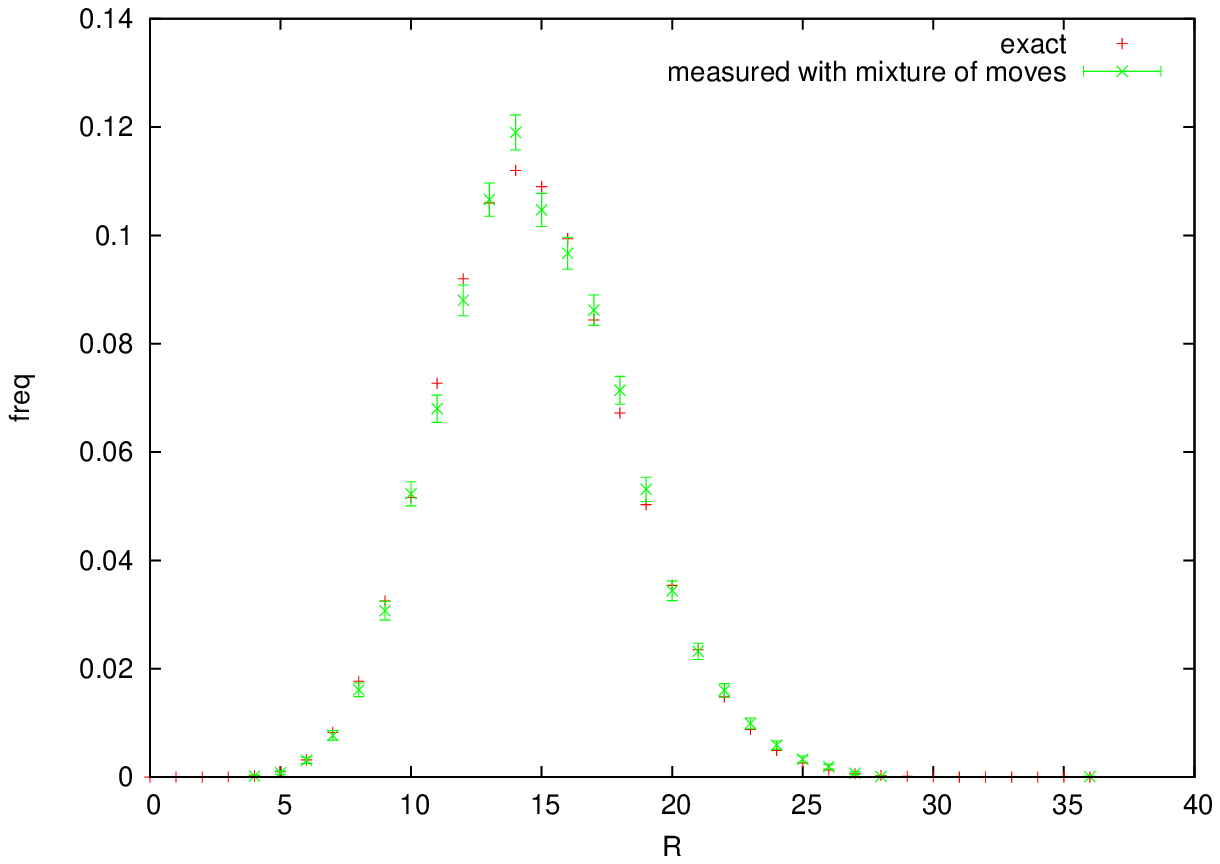}
  \end{minipage}
  \caption{\label{separate_moves.fig} Checking our code against numbers of
    relations for (naturally labeled) 9-orders.
    The green and blue data comes from 10k samples, the red is
    exact.  (Logarithmic scale on left, linear on right.)}
\end{figure}

Figure \ref{separate_moves.fig}
compares MCMC measurements of the number of relations
for 9-element posets against exact counts taken from the poset library.
The left graph separately tests the relation
and link moves, while the right tests a uniform mixture of both.  In each
case the measurements yield the correct values to within the estimated errors.
Notice that using 10k samples only allows us to measure frequencies down to
$10^{-4}$.
(The data-point with $R=36={9 \choose 2}$ represents the 9-chain that
was used as the starting configuration.)
\begin{figure}[htb!]
\centerline{\includegraphics{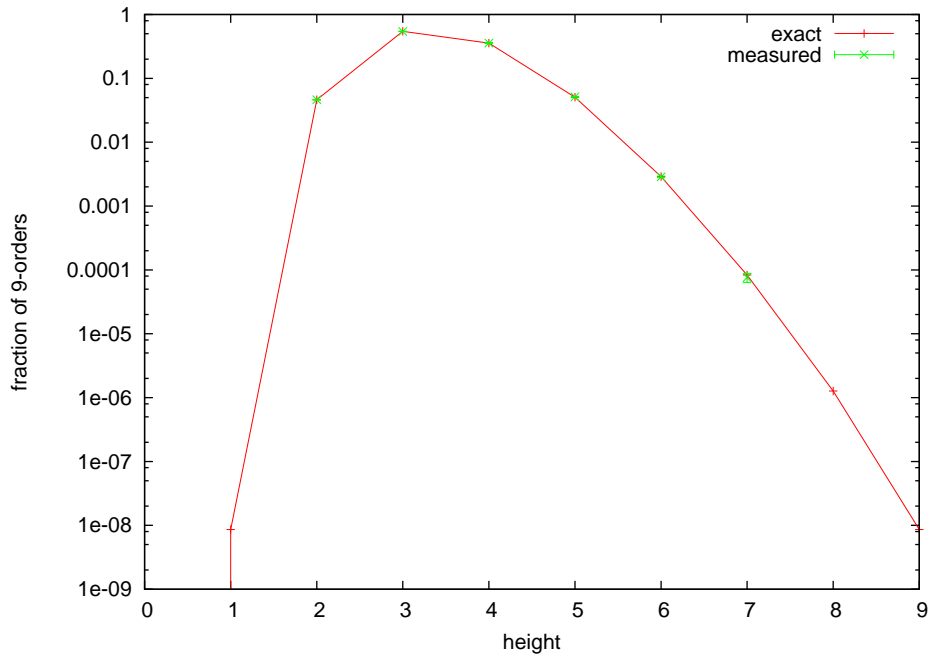}}
\caption{\label{height9works.fig} Heights of 9-orders  measured from 10k
    samples (green) compared with exact values (red). (Logarithmic scale on vertical axis.)}
\end{figure}
Figure \ref{height9works.fig} compares MCMC measurements of the heights of 9-element posets against
exact counts from the library.  Again the results match to within the error estimates.  Notice that
since there are fewer possible heights than relations, the measurements are more
accurate.  Of course, an MCMC measurement with 10k samples is still not able to resolve frequencies
$<10^{-4}$.

\section{Results}

How large must $n$ be for the asymptotic regime
described by Kleitman and Rothschild to set in?
As discussed in section \ref{s:onset},
we address this question
by measuring a number of order-invariants which
take on characteristic values in the KR regime, and which therefore
can herald the onset of this regime.
Among these are
the height,
the ordering fraction  $r$,
and the property of being ``layered'' 
(especially the ``connectedness'' expressed by condition (2) in section \ref{s:krtheorem}).  
We
also present results
on
the cardinalities of level 2
and of the minimal and maximal antichains.

We begin with a typical poset drawn from our
uniform sampling at
$n=82$, 
almost
the largest cardinality we have attempted to simulate.
Figure \ref{poset82.fig} is the Hasse
diagram of one such sample.
\begin{figure}
 \includegraphics[width=\textwidth]{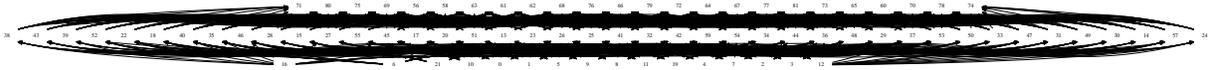}
 \caption{\label{poset82.fig} Sample 82 element poset.}
 \end{figure}
One sees that
the height is 3, and that the middle layer contains 42 elements,
exactly as one would expect for an asymptotic Kleitman-Rothschild poset.
By way of comparison, we show
in Figure \ref{poset20.fig}
a sample poset on 20 elements.
It clearly does not belong to the class
$\mathcal{A}(X_1, X_2, X_3)$
of the KR theorem of Section \ref{s:krtheorem}
(nor to $\mathcal{A}(X_1, X_2, X_3, X_4)$\footnote{$1 \prec 2 \prec 6 \prec 15$, so each would have to lie in the
corresponding layer $X_i$, $i=1\ldots 4$.  However $2 \nprec 18$, which
violates condition (2).}).

\begin{figure}[h!tbp]
 \includegraphics[width=\textwidth]{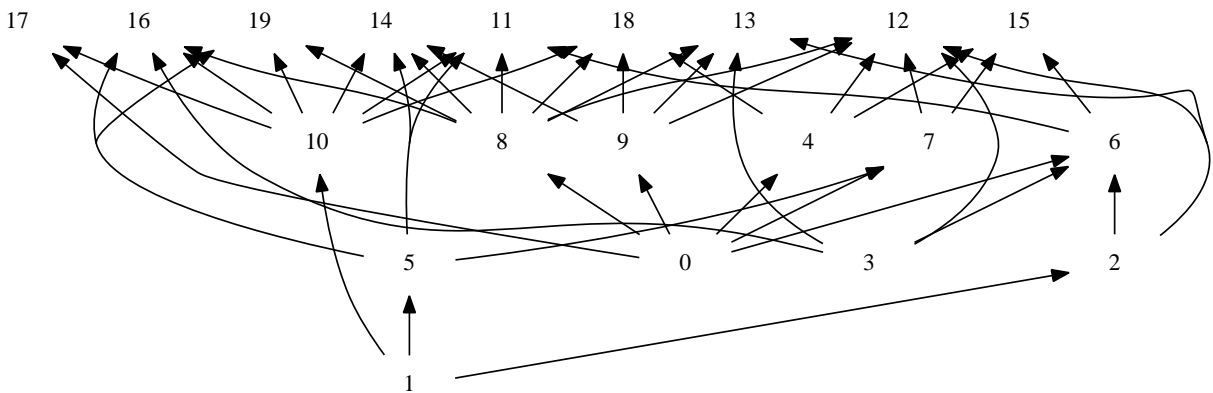}
 \caption{\label{poset20.fig} Sample 20 element poset.}
 \end{figure}

\subsection{Height}

As we have seen, the KR theorem implies that posets of height 3 dominate
as $n \rightarrow \infty$.
\begin{figure}[htb!]
\begin{minipage}[t]{.5\textwidth}
\includegraphics[width=\textwidth]{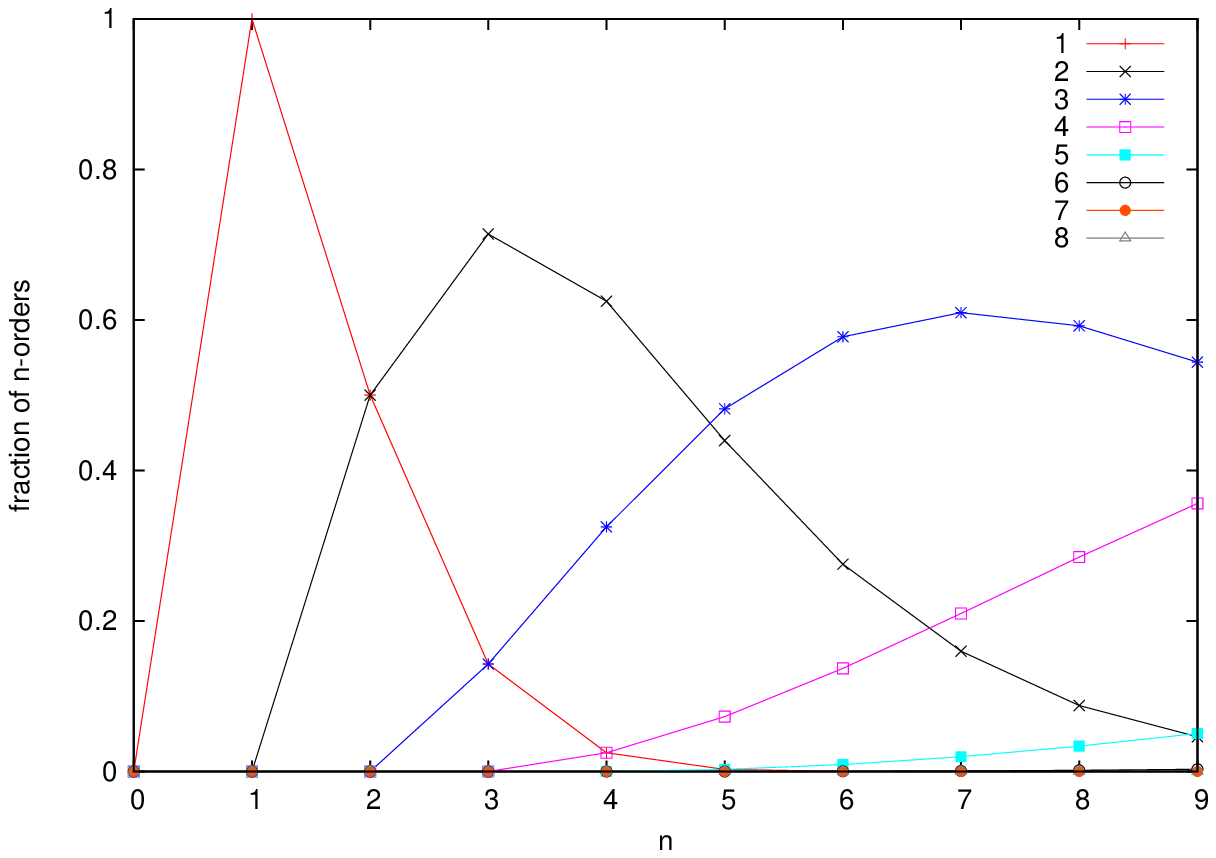}
\end{minipage}
\begin{minipage}[t]{.5\textwidth}
\includegraphics[width=\textwidth]{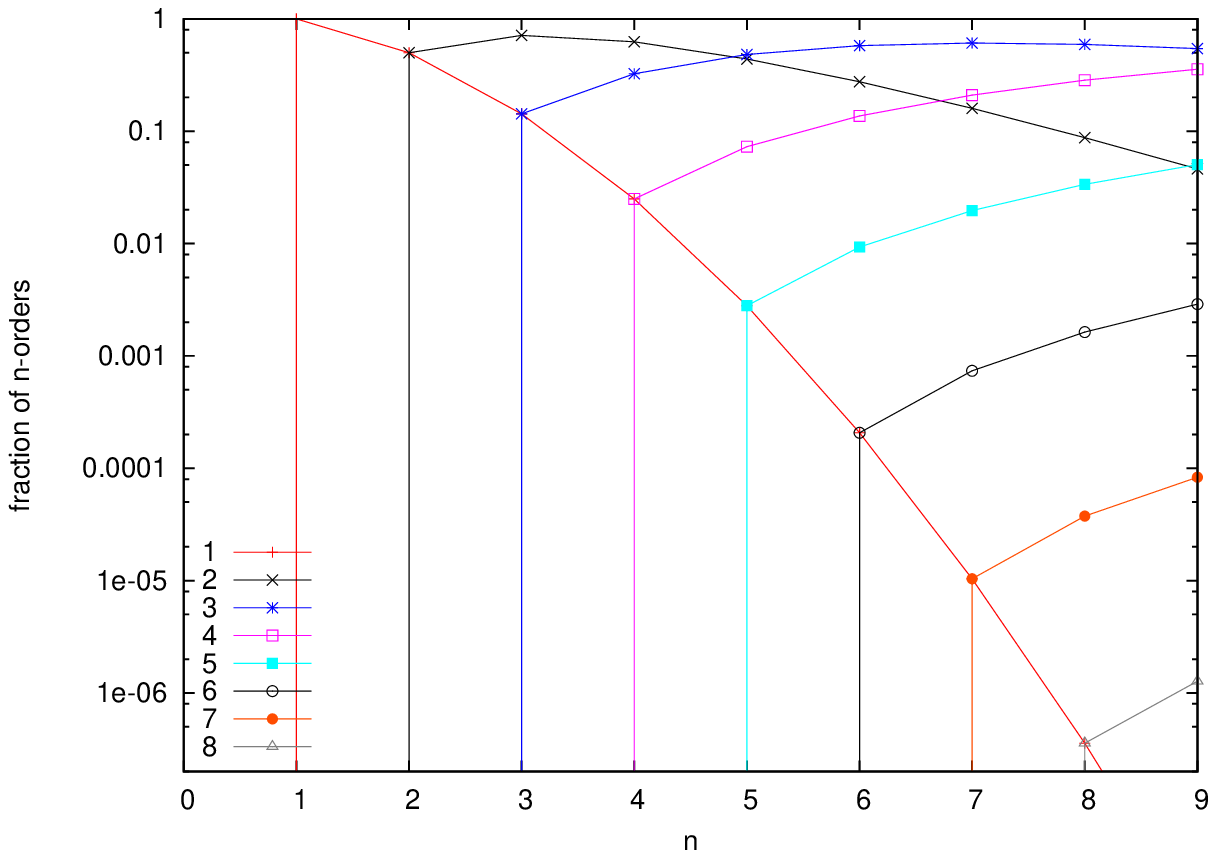}
\end{minipage}
\caption{\label{height9.fig} The exact fractions of $n$-orders with heights $h=1\ldots9$ as a
    function of $n$ for $n\leq 9\;$.
    (Linear scale on the left, logarithmic on the right.)}
\end{figure}
Figure \ref{height9.fig} portrays what we know
from exact computations \cite{smalln}
of the heights
of
posets with $n\leq 9$.
The fraction of
naturally labeled
posets with height $h=1 \ldots 9$
is plotted as a function of $n$.
One observes that while $n$-orders with height=3 are indeed in the majority for $n>5$,
their
relative
abundance begins to decrease at $n=7$,
while
the
 corresponding curves for heights $h\geq 4$ are all growing.
(This can be seen most clearly in the logarithmic plot on the right of figure \ref{height9.fig}.)
This indicates that, by any reasonable standard, $n=9$ is {\it\/not} in the asymptotic regime.

\begin{figure}[htbp!]
  \includegraphics[width=\textwidth]{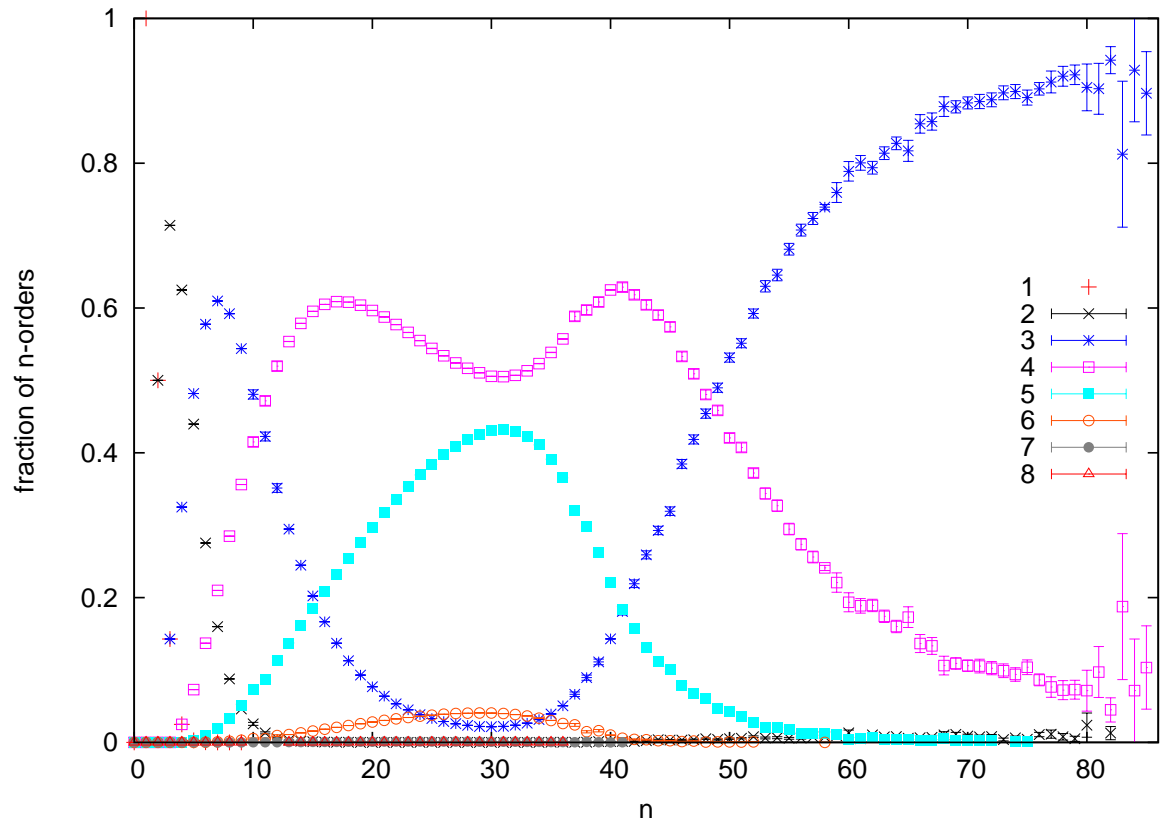}
  \caption{\label{heights.fig} Measured heights of posets with $n\leq 85$.}
\end{figure}

\begin{figure}[htbp!]
  \includegraphics[width=\textwidth]{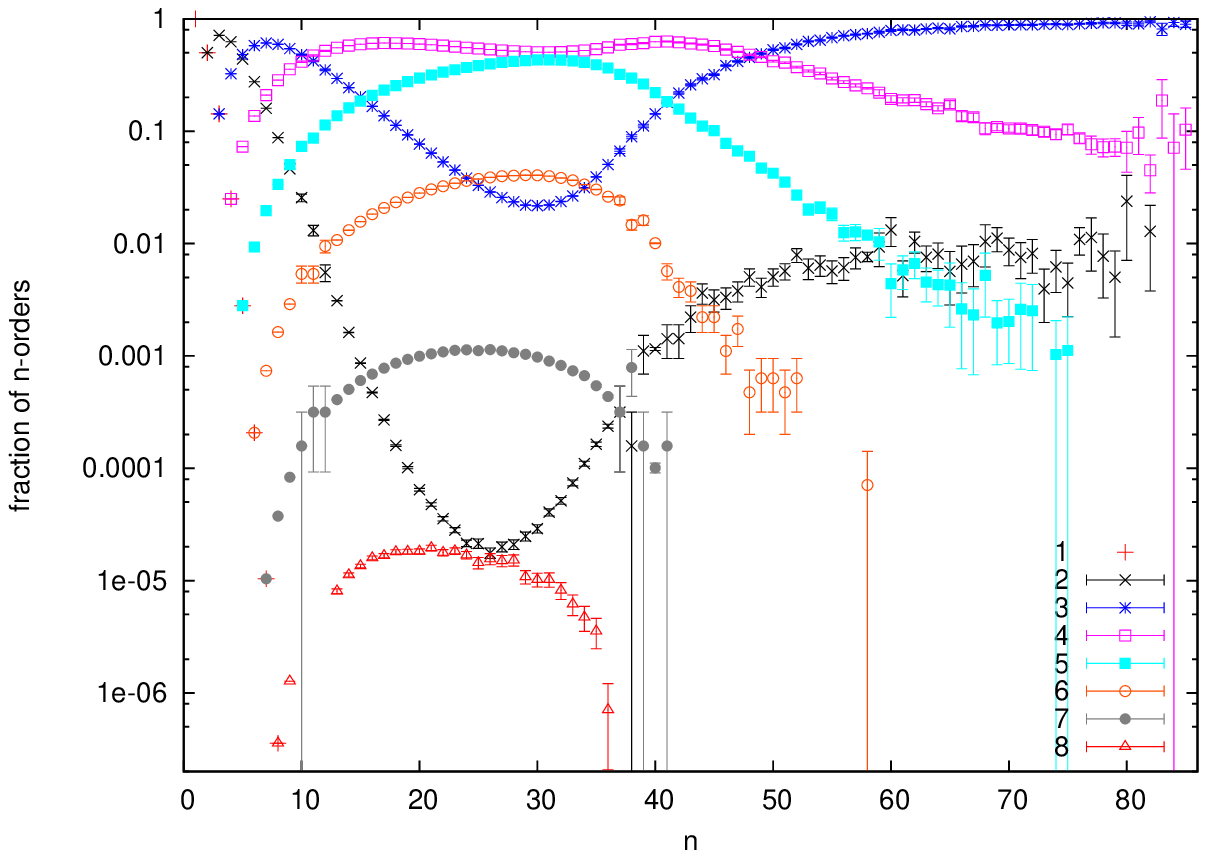}
  \caption{\label{heights-log.fig} Same as Figure \ref{heights.fig}, in logscale.}
\end{figure}
Our simulations now allow us to supplement these exact results with
statistical ones for $n\leq 85$.
We display our results in Figures \ref{heights.fig} and \ref{heights-log.fig},
which
extend
Figure \ref{height9.fig} out to $n=85$.
It is evident that the departure from KR-like behavior which started to
show up in the exact data continues for $n>9$ and
becomes most pronounced around
$n=30$.
The curves for $h\leq 3$ peak at small $n\leq 7$, and then diminish
rapidly as the $n\approx 30$ region is approached, whereas the height
fractions for $h\geq 5$ all reach their maxima at $n\approx30$,
strengthening the suggestion that something special occurs there.
Most remarkable are the curves for $h=3$ and $h=4$.
The latter exhibits two peaks,
one on each side of the $n\approx30$ region,
with a dip in the middle resembling those for heights two and three.
The $h=3$ curve, which asymptotically is supposed to rise to unity,
instead descends so precipitously after $n=7$ that for $25\leq n\leq 33$
the fraction of posets of height 3 falls below that of height 6,
decreasing to almost 1/50 at its minimum near $n=30$.
Not until $n\approx50$ does
$h=3$ meet and surpass $h=4$.

In the asymptotic regime,
all the height-fractions except that of height 3
should decay to zero,
with height 2 decaying the most slowly.
Our results are broadly
consistent with this expectation,
although they do not yet reveal any obvious sign of decrease for height 2.
\newpage 

In addition to plotting the full histogram of height as a function of $n$, it
is interesting to examine the mean height vs.\ $n$, as depicted in Figure
\ref{mean_height.fig}.

\begin{figure}[htbp]
  \includegraphics[width=\textwidth]{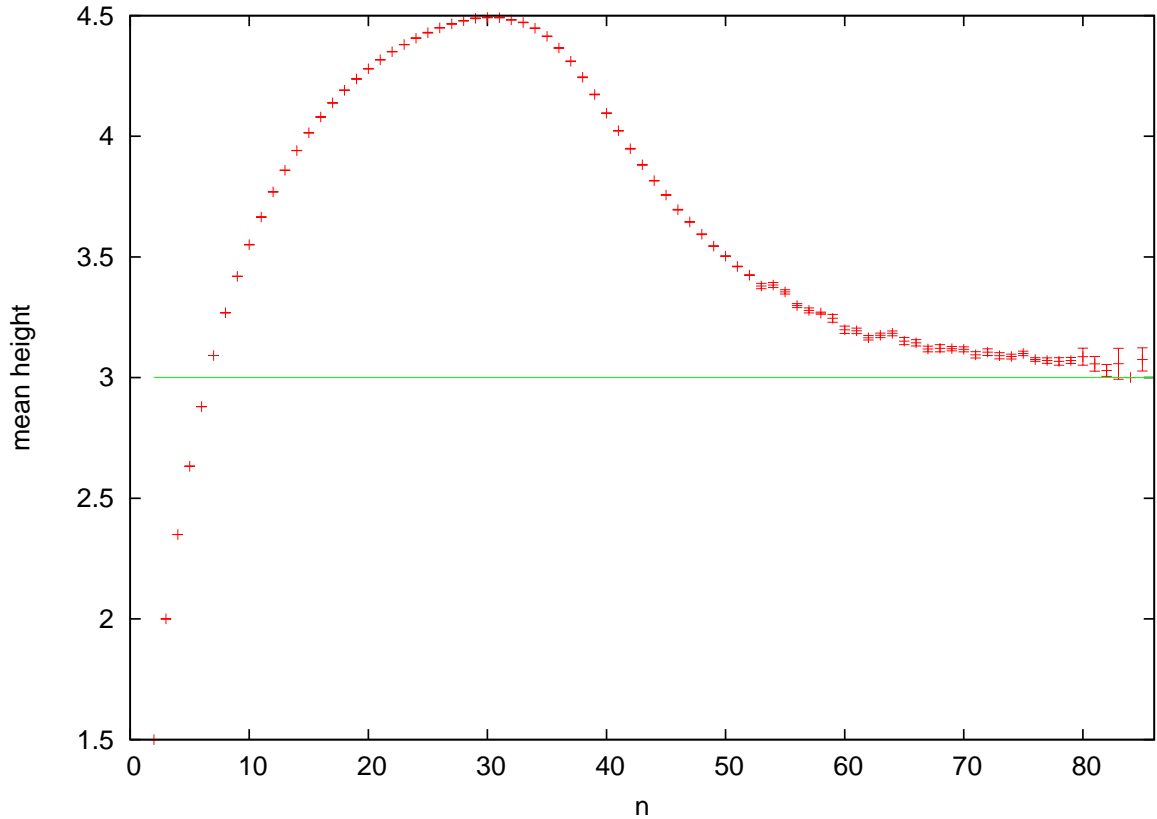}
  \caption{\label{mean_height.fig} Mean height as a function of $n$, compared
  against the asymptotic value 3.}
\end{figure}

\subsection{Connectivity of non-adjacent layers}

The KR theorem of section \ref{def.sec} informs us that asymptotically
almost every $n$-order belongs to the three-layer class of posets, and as
such fulfills conditions (1) and (2) given just above the theorem.
To examine the onset of this feature, we
define a characteristic function $\chi$ which takes the value 1 when conditions
(1) and (2) are satisfied for
the
partitioning of the poset into levels,\footnote{For computational convenience, we have used here levels as proxies for
  layers.  For this partitioning, condition (1) holds trivially.}
and is zero otherwise.
The fraction of such ``layered'' posets is depicted in Figure \ref{kr.fig}.
As with the relative abundances of $h=3$ and $h=4$, the layered posets
do not surpass the unlayered ones until about $n=55$.
\begin{figure}[htbp!]
  \begin{minipage}[t]{.5\textwidth}
  \includegraphics[width=\textwidth]{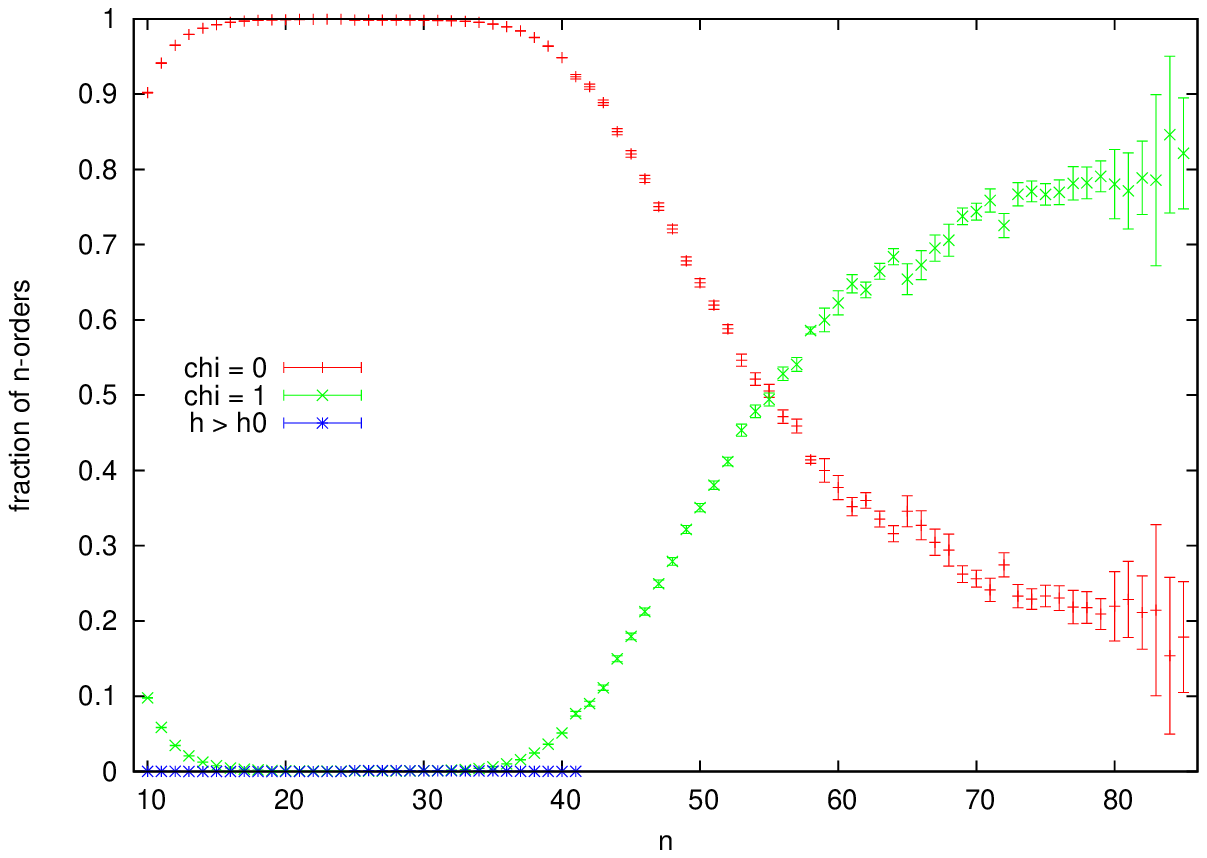}
  \end{minipage}
  \begin{minipage}[t]{.5\textwidth}
  \includegraphics[width=\textwidth]{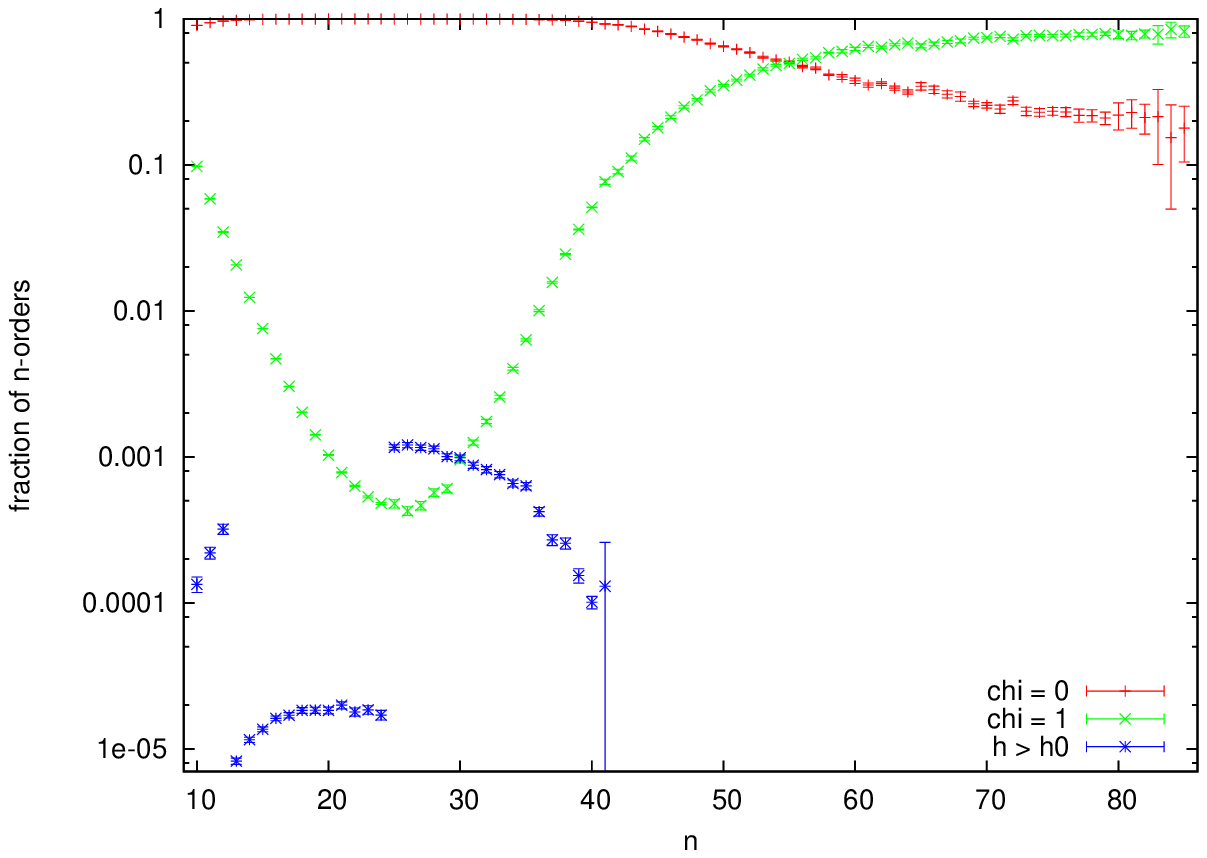}
  \end{minipage}
  \caption{\label{kr.fig} Fraction of posets
   whose levels do ($\chi=1$) and do not ($\chi=0$) satisfy
condition 
   (2)
   of Section \ref{def.sec}.
    Here, to speed the analysis,
    if the number of levels of
    the poset exceeds
    $h_0 =$ 7 (for $13\leq n \leq 24$) or 6 (otherwise) we abandon the
    check of that poset.  The plot on the right is in logscale.}
\end{figure}

A related condition, which is strictly stronger in the case of height 3
posets, is to demand that each of the minimal elements be related to each of
the maximal elements.  We do not check for this, since the
property of being layered is not yet satisfied by more than around 4/5 of the
85-orders.

\subsection{Cardinality of middle layer}

We now consider the cardinality of the middle layer.
Although it is not straightforward to determine whether a given poset
has a layered structure in the sense of Sec.\ \ref{def.sec}, it is easy
to compute its decomposition into levels.
Asymptotically, we expect level 2 to coincide with layer 2.
As an indication of how well the
quantitative bounds on layer size in the KR theorem are being satisfied,
we plot, in Fig.\ \ref{lev2.fig}, histograms of the cardinality of level 2
(divided by $n$) for $n=72$, 73, 76, and 79, restricted to the posets of
height 3.
As expected,
the curves
(which evidently differ little from each other)
are strongly peaked at $n/2$,
with a
width of less than $0.1$ and a
long tail that extends toward zero.
This is consistent with what one
would
expect from the Kleitman-Rothschild theorem for the cardinality of layer 2.
\begin{figure}
\centerline{\includegraphics{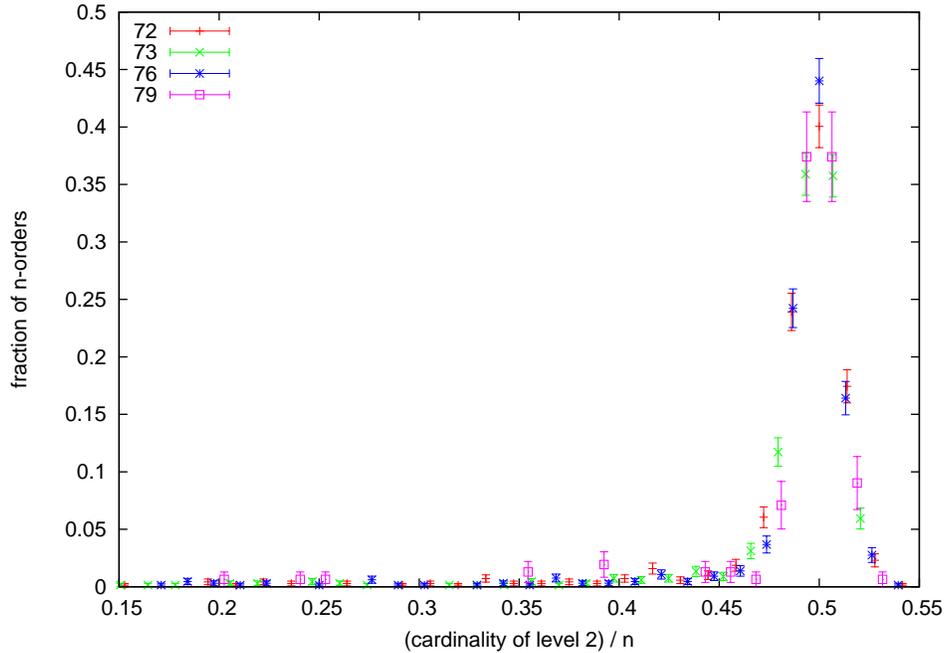}}
\caption{\label{lev2.fig} Histograms for the cardinality of level 2 for $n=$
  72, 73, 76, and 79, restricted to the posets of height 3.
}
\end{figure}

\subsection{Time asymmetry}
\label{time_asym.sec}
Define the time reverse of a poset $P$ as a poset $P'$ which is identical to $P$ except that the
direction of each relation is reversed: $x\prec y$ in $P'$ iff $y\prec x$ in $P$.
The version of the KR theorem we have quoted above, implies that a
typical $n$-order will be almost perfectly time-symmetric for $n\gg1$.
However, we have argued that our weighting by natural labellings is
liable to invalidate this sort of conclusion.  Indeed, an unconstrained
division of elements between the the top and bottom layers would imply
a significant amount of asymmetry, even as $n\to\infty$.
In Fig.\ \ref{asym58.fig} we examine the situation for $n=58$, by measuring a
histogram of the number of minimal elements $|\mathrm{min}|$,
the number of maximal elements $|\mathrm{max}|$,
and the difference $|\mathrm{max}| - |\mathrm{min}|$.
As they must, the histograms of $|\mathrm{min}|$ and $|\mathrm{max}|$
match each other (to within the estimated errors).
But we notice that the cardinalities of the minimal and maximal layers
are more likely to differ by $\approx 17$ than by smaller or bigger values.
In this sense most naturally labeled 58-orders exhibit a substantial time
asymmetry, exceeding what one might have expected from the naive KR
bounds.
The same sort of asymmetry persists for other values of $n$ as well.

\begin{figure}
\psfrag{|min|}{\footnotesize $N_\mathrm{min}$}
\psfrag{|max|}{\footnotesize $N_\mathrm{max}$}
\psfrag{|max| - |min|}{\footnotesize $N_\mathrm{max} - N_\mathrm{min}\;\;\;\;$}
\centerline{\includegraphics{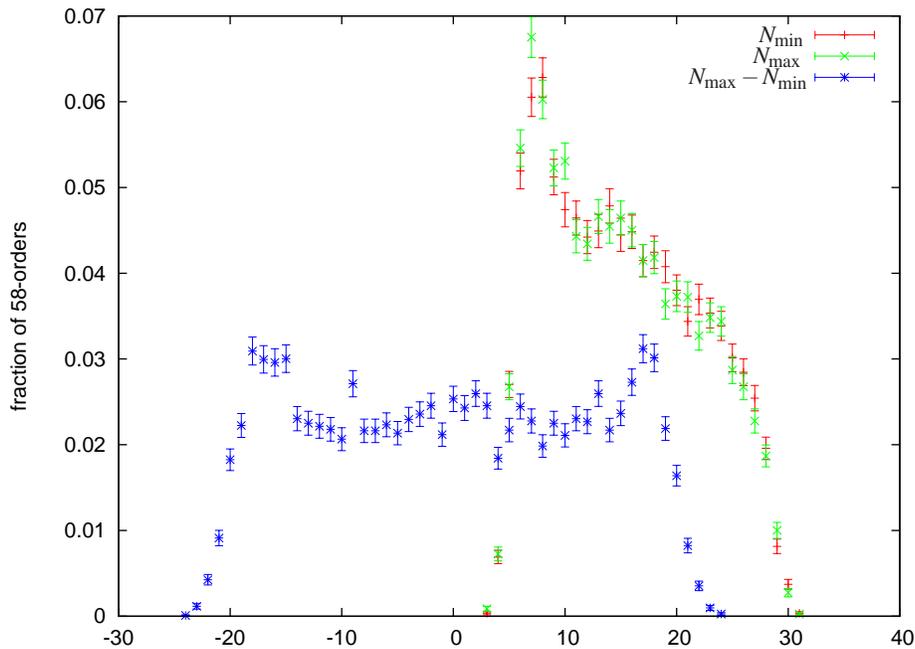}}
\caption{\label{asym58.fig} Histograms of the number of minimal and maximal
  elements, and their difference, for posets of 58 elements.}
\end{figure}

Interestingly, a similar
asymmetry was found in a generalization of
the uniform
model of random partial orders in
which the ordering fraction $r$
is held fixed~\cite{dhar}.

\subsection{Ordering fraction}

The ordering fraction $r$ of a poset is by definition the number $R$ of
pairs of related elements $x \prec y$ divided by ${n \choose 2}$.
As described earlier in connection with the KR theorem, we expect --- for
naturally labeled orders --- that $r$ asymptotically will be distributed
between 1/4 and 3/8 with a square-root-divergent peak at $r=3/8$.

Figure\ \ref{r.fig} presents a histogram of $r$ for $n=58$.
Evidently, it is nicely consistent with these expectations.
The peak at $r=3/8=0.375$ is present,
and the cutoff at $r=1/4$ is also observed.
Except for the fact that the histogram extends some way to the right of
the peak, the agreement is
surprisingly good,
considering that a large number of posets of height 4 are still present
at $n=58$, as we saw above.
\begin{figure}
  \centerline{\includegraphics{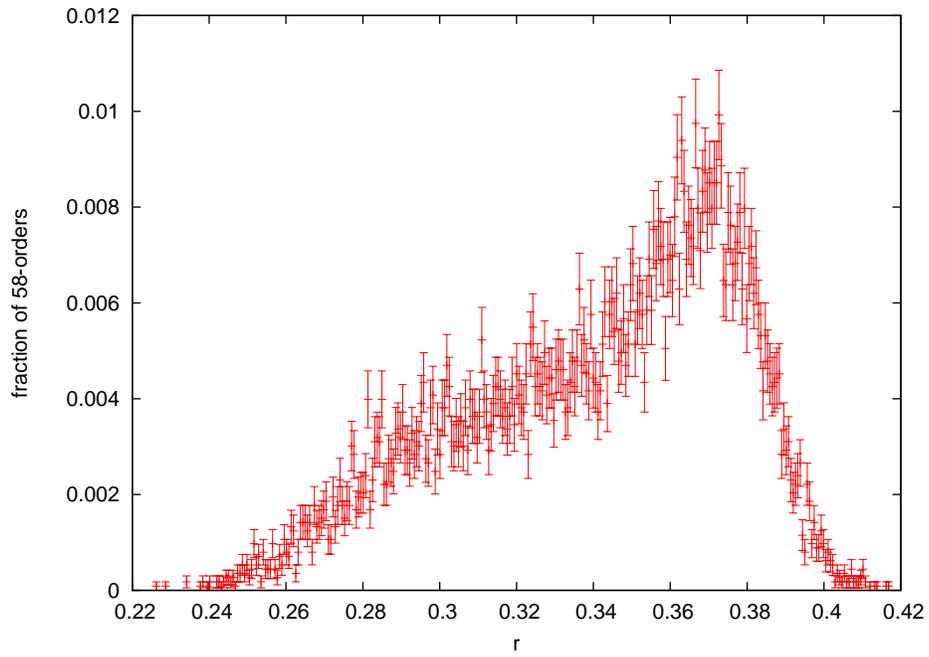}}
  \caption{\label{r.fig} Histogram of ordering fraction $r$ for $n=58$.}
\end{figure}

After 3-layer orders (or as a special case thereof), 2-layer or
``bipartite'' orders are expected to be next most abundant for large
$n$.  A typical bipartite order (one with both layers of size
$\approx{n/2}$, such that each element of one layer is related to
approximately $n/4$ elements of the opposite layer) has $r\approx1/4$,
furnishing in effect the lower end of the 3-layer histogram.
Once again, this is quite compatible with what one observes
in Figure~\ref{r.fig} in the neighborhood of $r=0.25$.

\begin{figure}
  \centerline{\includegraphics{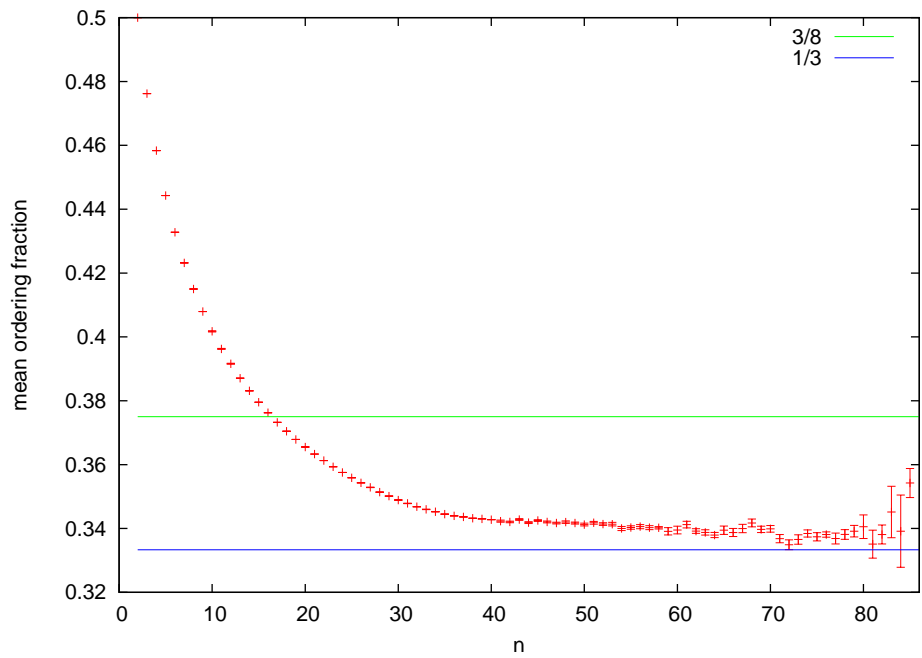}}
  \caption{\label{r_mean.fig} Plot of mean ordering fraction $\langle r \rangle$
  versus $n$.}
\end{figure}
To get additional feel for how $r$ varies with $n$, we computed its
expectation value, as shown in Figure \ref{r_mean.fig}.
It seems evident that the mean is converging to something very close to
the expected\footnote{``expected'' but not ``predicted'', since we understood the significance of
 natural vs. arbitrary labellings only after seeing this curve!}
value of 1/3, and not to the ``unlabeled KR'' value of 3/8.

\section{Conclusions}

Thanks to the groundbreaking results of Kleitman and Rothschild
the asymptotic
structure and
enumeration of
the finite partial orders
are fairly well understood.
However,  almost 40 years since
the publication of their main results,
relatively little is known about
the onset of the asymptotic
regime,
i.e., the poset size at which their results become relevant.
While small posets
can be generated
exhaustively
on a computer,
the number of $n$-orders grows so rapidly with $n$ that one cannot
access the asymptotic, KR regime in this way.

The Markov chain algorithm presented herein
opens
a new window
onto
this question
by sampling uniformly from the set of naturally labeled partial orders.
We have used it to explore the structure of a ``typical'' poset up to
$n=85$.
Although
we still do not reach the deep asymptotic regime in this way, our results
are consistent with a monotonic approach to
asymptopia
 which begins
around $n=45$ and is in full swing by $n=80$.
It seems likely from our
data
that
this trend
simply continues until it reaches substantial completion well
past $n=100$.
Our simulations also reveal the presence of an interesting intermediate
regime for $n<50$, which is quite different from the asymptotic one, and
whose best characterization remains to be understood.

Even though our algorithm allows uniform sampling of partial orders which are much larger than
have been considered previously, one would like to be able to
simulate still larger
posets
which are deep within the asymptotic regime.
The construction of efficient algorithms, however, is
hampered by the constraint of transitivity,
which entails complicated dependencies among the defining relations $x\prec y$.
One could circumvent this difficulty by identifying a set of
simple,
independent
variables in terms of which the partial order could be defined, however we are
not aware of any such representation of (generic) partial orders.
To bypass
the problems arising from the transitivity constraint,
one might try to devise
`cluster algorithms' analogous to the Wolff algorithm for spin
systems~\cite{wolff},
which could
change a large number of relations in a single move,
while still maintaining transitivity.
An alternative approach would be to drop transitivity and
expand the
sample space
from $\Omega_n$ to the set of all directed acyclic
graphs.
One could then
impose Boltzmann weights
which would keep the violations of transitivity small,
or alternatively
which would compensate for the overcounting resulting from failing to
take a transitive closure.\footnote{The observables would still be measured on the transitive closure of
  the directed graph, its ``projection down to $\Omega_n$''.}
The latter approach is currently being pursed \cite{roysurya}.
Preliminary results for $n<55$ are in good
agreement with the results reported above.
In particular, the dip in the
profiles for $h\leq4$
appearing
in Figures\ref{heights.fig} and \ref{heights-log.fig} seems to be
reproduced
with high precision.

Throughout we have considered the set of naturally labeled partial orders.
This choice was largely for convenience.  Labeled posets are
conveniently represented on a
computer by their
adjacency matrices,
without any need to search for
isomorphisms between different labellings.
The choice is also convenient because
one can easily take advantage of the CausalSets Toolkit within the Cactus HPC
framework, which
currently
assumes natural labellings throughout.

We conclude with
some
questions
concerning
the possible relevance of posets
to the still unsolved mystery surrounding the so-called arrow of time.
If all the known laws of physics are invariant under reversal of the direction of time,
from whence comes the extreme time-asymmetry that we observe?
Locally finite posets
in which the partial ordering is interpreted temporally
(and which are called in that context causal sets)
have been proposed as
the deeper physical structure from which spacetime
geometry emerges,
and with discreteness comes a well-defined ``counting entropy''.
Could it be this entropy, rather than the kind of low entropy initial
conditions on matter, radiation and event horizons, that are usually
invoked, which will account for the arrow of time?
Could time-asymmetric universes simply be much more plentiful than time
symmetric ones, because time-asymmetric posets are
more plentiful?
A hint of this comes from the asymmetry observed in Sec.\ \ref{time_asym.sec},
but it's doubtful that simple counting could by itself suffice to explain
the observed degree of asymmetry in the cosmos.
But perhaps in conjunction with other dynamical constraints,
like those in \cite{dhar},
it could be a part of the answer.

\paragraph{Acknowledgments}

We are grateful to a number of people for comments and suggestions, including
Samo Jordan, Lisa Glaser, Andrzej Gorlich, Denjoe O'Connor, Jeff Remmel, Orest Bucicovschi,
and Graham Brightwell.

This work was funded by a grant from the Foundational Questions Institute
(FQXi) Fund  on the basis of proposal FQXi-RFP3-1018.  This work was also
supported in part under an agreement with Theiss Research and funded by a
grant from the FQXI Fund on the basis of proposal FQXi-RFP3-1346 
to the Foundational Questions Institute.
The FQXI  Fund is a donor advised fund of the Silicon Valley Community
Foundation.
This research was also supported in part by NSERC through grant RGPIN-418709-2012.
This research was additionally supported in part by Perimeter Institute for Theoretical Physics.
Research at Perimeter Institute is supported by the Government of Canada
through Industry Canada and by the Province of Ontario through the
Ministry of Research and Innovation.
JH also receives support form EPSRC grant \textit{DIQIP} and ERC grant \textit{NLST}.
This material is based in part upon work supported by DARPA under
Award No.\ N66001-15-1-4064.  Any opinions, findings, and conclusions or
recommendations expressed in this publication are those of the authors
and do not necessarily reflect the views of DARPA.
This work used the Extreme Science and Engineering Discovery Environment (XSEDE), which is supported
by National Science Foundation grant number ACI-1053575.

We thank Yaakoub El Khamra for his
encouragement and assistance with profiling our code on Lonestar
(\url{www.tacc.utexas.edu/resources/hpc/lonestar}).
Early trials were conducted on the HPC cluster at the Raman Research Institute.

\appendix
\section{Appendix}
We here crudely estimate the number of naturally labeled 3-layer orders
with layer sizes $n_1$, $n_2$, $n_3$.

Assume
(as will almost always be true)
that the poset is asymmetric (automorphism free), and for convenience
let its elements be distinguishable.
The number of ways to fill in the links between adjacent layers is
$$
  2^{n_1 n_2 + n_2 n_3} = 2^{(n_1+n_3)n_2} = 2^{n_2 (n - n_2)} \;.
$$
However, permuting the elements within any single layer produces
a different set of links, but an isomorphic poset.  Compensating
for this over counting yields, for unlabeled 
posets,
$$
  {2^{n_2(n-n_2)} \over  n_1! \; n_2! \; n_3! } \;.
$$

Now let $L(a,b)$ be the average number of natural labellings of a
bipartite order with layers of size $a$ and $b$.
From \cite{brightwell2} comes the estimate
$$
      L(a,b) = c \ a! \ b! \;,
$$
where $c=c(a,b)$ depends on $a$ and $b$, but goes over to the constant
$\eta=3.4627\dots$ as $n\to\infty$  (and is typically smaller than that
for finite $n$).  Thus we estimate
$$
      c(n_1,n_2) \  c(n_2,n_3) \  n_2! \  2^{n_2(n-n_2)}
$$
naturally labeled orders with layer-sizes $n_1, n_2, n_3$, where
asymptotically $c=\eta$.

The exponential factor peaks sharply at a middle-layer size of
$n_2=n/2$, and asymptotically dominates everything else.  What remains
is then independent of $n_1$ and $n_3$, to the extent that $c$ and
$n_2!$ can be treated as constant.

\input{KROnset5.1.bbl}
\end{document}

(prog1 'now-outlining
 (Outline* 
     "\f"                  1
     "     "     "     "\\documentclass"     1
     "\\begin{doc"         1
     "\\end{doc"           1
     "\\title"             1 
     "\\author"            2 
     "\\date"              2  
     "\\abstract"          1
     "\\Abstract"          1
     "\\begin{abstract"    1
     "\\section"           1
     "\\subsection"        2
     "\\appendix"          1
     "\\ReferencesBegin"   1
     "     "\\ref "              2
     "\\begin{thebibliog"  1
     "\\bibitem"           2
   ))

================================================================================

%% file: mydefs.tex
\def\be{\begin{displaymath}}
\def\ee{\end{displaymath}}
\def\bne{\begin{equation}}
\def\ene{\end{equation}}
\def\bee{\begin{eqnarray*}}
\def\eee{\end{eqnarray*}}
\def\bnee{\begin{eqnarray}}
\def\enee{\end{eqnarray}}

\font\openface=msbm10 at10pt

\def\past{\mathrm{past}}
\def\ipast{\mathrm{incpast}}
\def\fut{\mathrm{fut}}
\def\ifut{\mathrm{incfut}}

\def\union{\cup}
\def\nc2{{N \choose 2}}
\def\Nc2{{N \choose 2}} 

\def\iff{\Leftrightarrow}

\def\Pr{\mathrm{Pr}}

\def\lto{\mathop
        {\hbox{${\lower3.8pt\hbox{$<$}}\atop{\raise0.2pt\hbox{$\sim$}}$}}}
\def\gto{\mathop
        {\hbox{${\lower3.8pt\hbox{$>$}}\atop{\raise0.2pt\hbox{$\sim$}}$}}}



%% file: journals.tex
\newcommand{\eprint}[1]{$\langle$e-print arXiv: #1$\rangle$}